\PassOptionsToPackage{dvipsnames}{xcolor}
\documentclass[nohypdvips,11pt]{siamart220329}

\usepackage{enumitem}
\usepackage{graphicx,epstopdf} 
\usepackage[caption=false]{subfig} 
\usepackage{multicol,multirow}
\usepackage{longtable}
\usepackage{lscape}
\usepackage{float}
\usepackage{physics}
\usepackage{lipsum}
\usepackage{amsfonts}

\usepackage{tikz}
\usepackage{booktabs}
\usepackage[algo2e,boxruled,lined,linesnumbered]{algorithm2e}
\usepackage[numbers]{natbib}

\newsiamremark{remark}{Remark}
\newsiamremark{hypothesis}{Hypothesis}
\crefname{hypothesis}{Hypothesis}{Hypotheses}
\newsiamthm{claim}{Claim}

\title{DRSOM: A Dimension Reduced Second-Order Method \thanks{Submitted to the editors DATE.
        \funding{This research is partially supported by the National Natural Science Foundation of China (NSFC) [Grant NSFC-72150001, 72225009, 11831002]
        }}}

\author{
    Chuwen Zhang\thanks{Shanghai University of Finance and Economics,
        (\email{chuwzhang@gmail.com}).}
    \and Dongdong Ge\thanks{Corresponding author. Shanghai University of Finance and Economics, \\
        \indent(\email{ge.dongdong@mail.shufe.edu.cn}), }
    \and Chang He\footnotemark[2]
    \and Bo Jiang\footnotemark[2]
    \and Yuntian Jiang\footnotemark[2]
    \and Yinyu Ye\thanks{Stanford University,
        (\email{yinyu-ye@stanford.edu}).}
}

\newcommand{\mn}{\mathsf{\min}~}
\newcommand{\st}{\mathrm{s.t.}~}

\newcommand{\rea}{\mathbb{R}}

\newcommand{\xk}{x_k}
\newcommand{\Hk}{H_k}
\newcommand{\Qk}{Q_k}
\newcommand{\gk}{g_k}
\newcommand{\dk}{d_k}
\newcommand{\gkn}{g_{k+1}}
\newcommand{\xkn}{x_{k+1}}
\newcommand{\dkn}{d_{k+1}}
\newcommand{\lk}{\lambda_k}

\newcommand{\hkap}{\tilde{H}_k}
\newcommand{\vk}{V_k}
\newcommand{\spa}{\mathsf{span}}
\newcommand{\spak}{\mathcal L_k}
\newcommand{\ak}{\alpha_k}

\newtheorem{thm}{Theorem}
\newtheorem{assm}{Assumption}

\newtheorem{lem}{Lemma}

\newcommand{\drsom}{\textrm{DRSOM}}

\newcommand{\lbfgs}{\textrm{LBFGS}}
\newcommand{\newtontr}{\textrm{Newton-TR}}
\newcommand{\cg}{\textrm{CG}}
\newcommand{\gd}{\textrm{GD}}

\begin{document}

\maketitle

\begin{abstract}%
  In this paper, we propose a Dimension-Reduced Second-Order Method (DRSOM) for convex and nonconvex (unconstrained) optimization. Under a trust-region-like framework, our method preserves the convergence of the second-order method while using only curvature information in a few directions. Consequently,
  the computational overhead of our method remains comparable to the first-order such as the gradient descent method.
  Theoretically,
  we show that the method has a local quadratic convergence and a global convergence rate of \(O(\epsilon^{-3/2})\) to satisfy the first-order and second-order conditions if the subspace satisfies a commonly adopted approximated Hessian assumption.
  We further show that this assumption can be removed if we perform a \emph{corrector step} using a Krylov-like method periodically at the end stage of the algorithm.
  The applicability and performance of DRSOM are exhibited by various computational experiments, including $\mathcal L_2 -\mathcal L_p$ minimization, CUTEst problems, and sensor network localization.
\end{abstract}



\section{Introduction}
In this paper, we consider the following unconstrained optimization problem
\begin{equation}\label{Prob:main}
  \min_{x\in \rea^n} f(x),
\end{equation}
where \(f:\rea^n \to \rea\) is
twice differentiable and possibly nonconvex and \( f_{\inf}:=\inf f(x) > -\infty\).
We aim to find a  ``stationary point'' $x$ such that
\begin{equation}\label{cond.eps.grad}
  \|\nabla f(x)\|\le \epsilon
\end{equation}
and $x$ approximately satisfies the second-order necessary conditions in a certain subspace.
Historically speaking, various methods have been proposed to solve \eqref{Prob:main}, among which the first-order methods (FOM) based on the gradient of $f$ are the most popular ones.
A first-order method computes a descent direction that iteratively decreases the function value, in which the negative gradient, the momentum \cite{polyak_methods_1964} can be applied as to construct the iterates.
Theoretically, FOMs has the global complexity rate in \(O(\epsilon^{-2})\) \cite{nesterov_lectures_2018}, however, they only converge to a point that satisfies the first-order condition \eqref{cond.eps.grad}.


To escape the first-order saddle point, people may resort to the second-order methods based on the Hessian of $f$ with variants of Newton-type iterations. Traditionally, the global convergence of second-order method can be ensured by frameworks like line search and trust-region methods \cite{nocedal_numerical_2006}. While these classical methods are efficient in practice due to local superlinear convergence rate, the global iteration complexity of the standard Newton's method is generally in the same order as that of the gradient descent method \cite{cartis_complexity_2010}.
To clear the mist, 
\citet{nesterov_cubic_2006} established the first \(O(\epsilon^{-3/2})\) iteration complexity for second-order methods using cubic-regularized Newton method (also see \citet{cartis_adaptive_2011,cartis_adaptive_2011-1}). However, the complexity analysis for the trust-region methods seems to be more challenging. \citet{ye_second_2005} in his early lecture notes showed that the trust-region method with a fixed-radius strategy also has the iteration complexity of \(O(\epsilon^{-3/2})\), which matches the iteration bound of the aforementioned cubic regularized Newton method. A more practical \emph{adaptive} trust-region method with $O(\epsilon^{-3/2})$ bound was only unveiled in \cite{curtis_trust_2017} until very recently.

When the problem dimension gets larger, both explicitly computing the Hessian matrix and solving Newton-type subproblems could be more costly. A practical implementation of the second-order method usually uses a Lanczos method to iteratively find inexact solutions \cite{conn_trust_2000,cartis_adaptive_2011-1,curtis_worst-case_2022}, where the ``inexactness'' can be controlled by the quality of Hessian approximation in some sense \cite{dennis_quasi-newton_1977, cartis_adaptive_2011, curtis_trust_2017,xu_newton-type_2020}. However, even this approach can be impractical for huge-sized problems, e.g., solving sensor-network localization for thousands of sensors that may induce a undirected graph with millions of edges.
By contrast, the first-order methods that only require gradient and cheap directions are scalable well to those huge-sized problems. However, their numerical performance heavily relies on \emph{proper} step sizes and other hyperparameters, e.g., using various line search methods \cite{hager_algorithm_2006}, and adaptive strategies \cite{kingma_adam_2014} in stochastic settings. The choice of step size itself requires plentiful experiment and trial procedure.  Despite all these difficulties, these methods are prevalent in large scale optimization.

Our approach in this paper is motivated from finding stepsizes for adaptive gradient methods using the Polyak's momentum \cite{ye_second_2005}. Specifically, we introduce a \emph{Dimension-Reduced Second-Order Method} (DRSOM)
to construct and solve a low-dimensional quadratic model to select the stepsizes.
The extra overhead of DRSOM is mostly due to constructing and solving such quadratic subproblems.

\subsection{Related Work}

DRSOM adopts the basic idea of manipulating the low-dimensional subspace. In optimization, the subspace method is fundamental and widely studied.  For example, the gradient method can be seen as performing a one-dimensional subspace minimization. In trust-region literature, such a minimizer in $\spak = \spa\{\gk\}$ is referred to as the Cauchy point \cite{conn_trust_2000}. In comparison,
\citet{yuan_subspace_1995} discussed the quadratic approximation over different simple and specific subspaces with an emphasis on their close relationship to the (nonlinear) conjugate gradient method. The authors also provided convergence proofs for the cases based on exact and inexact linear searches.
A notable extension of the idea in \cite{yuan_subspace_1995} was made in \cite{wang_subspace_2006} to the trust-region methods. Briefly speaking, if the subspace is constructed from the previous gradients under the quasi-Newton framework, then one can show that the subspace minimization generates exactly the same iterates of a full-dimensional quasi-Newton algorithm (see \cite[Lemma 2.3]{wang_subspace_2006}). A comprehensive summary on subspace methods can be found in \cite{yuan_review_2014, liu_subspace_2021}.

In contrast with the deterministic setting mentioned above, the \emph{sketching techniques} can also be seen as a \emph{stochastic} subspace method. Early results of sketching appear in solving least-square problems and more generally in numerical linear algebra \cite{woodruff_sketching_2014}. The basic idea is to replace the full Hessian with its approximation     ---without destroying the convergence. Using different sketching schemes on the Hessian matrix, \citet{pilanci_newton_2017} discussed the convergence rates of the Newton method. This idea later became popular in stochastic second-order optimization \cite{lacotte_adaptive_2021,berahas_investigation_2020,liu_subspace_2021}, and in natural gradient method \cite{yang_ng_2021, yang_sketch-based_2022} for training deep neural networks. In essence, since the sketching methods use random matrices, the convergence and complexity analyses are typically probabilistic or as expectation in some error metrics. Moreover, the assumptions required by sketching methods, once again, are usually expressed in expectations or probabilistic. In a high-level view, the random sketch matrix has many notable advantages, for example, it may help the algorithm to escape from saddle points efficiently in nonconvex optimization. In comparison, DRSOM is greedier and straightforward since the subspace is kept at a minimum behavior and only increased when using the corrector steps.

As for the momentum-based gradient method, its global convergence rate was only established for convex optimization problems, in terms of the average iterates \cite{ghadimi_global_2015} as well as the last iterate \cite{sun_non-ergodic_2019}. When the objective function is nonconvex, the sequence generated by the heavy ball method converges to a first-order stationary point \cite{zavriev1993heavy, ochs2014ipiano}. Later, \citet{sun2019heavy} proved that heavy ball method could even escape saddle point. However, the complexity rates are generally not superior to the gradient method under the nonconvex setting \footnote{To our best knowledge, we are unaware of better complexity results of Polyak momentum in nonconvex settings}. A recent work in \cite{apidopoulos_convergence_2022} studied the convergence of heavy ball method under the Polyak–Łojasiewicz condition based on an analysis of the dynamic system.

An earlier line of \emph{modified first-order methods} in the nonconvex case explored the idea of using second-order information in the Nesterov's accelerated gradient method (AGD), where the momentum direction is also appeared. For example, \citet{carmon_accelerated_2018,carmon_convex_2017} used the negative curvature with AGD to achieve a better complexity bound of \(O(\epsilon^{-7/4}\log(\epsilon^{-1}))\) for first-order methods. In this vein of research, the improved bound is based on finding the negative curvature with less dimension-dependent complexity\footnote{For example, one can use a randomized Lanczos method which has a complexity bound of \(O(\epsilon^{-1/4}\log(\epsilon^{-1}))\)}. With the help of negative curvature, these methods are able to ``escape'' from saddle points efficiently and converge to the approximate second-order points in sharp comparison to classical first-order methods.

\subsection{Our Contribution}
In light of the above discussion,
this paper extends previous results by demonstrating that one may efficiently enjoy the \emph{charisma} of second-order information while keeping the low per-iteration cost of a first-order method.

\begin{itemize}[leftmargin=*]
  \item Firstly, we introduce a general dimension-reduced framework DRSOM by solving cheap subproblems without
        invoking the full Hessian matrix, see \autoref{alg.concept}. We further propose some efficient constructions of the low dimensional subproblems based on Hessian-vector products or \emph{interpolation}.
        Notably, the DRSOM can be viewed from the subspace perspective \cite{yuan_subspace_1995,yuan_recent_2015,liu_subspace_2021}, however the associated analysis of global and local rates of convergence is missing.
        To our best knowledge, the interpolation technique is rarely celebrated in subspace and dimension reduction literature.
  \item Secondly, we analyze the convergence of DRSOM based on the discussion in \cite{ye_second_2005}. Under a commonly adopted approximated Hessian assumption (cf. \cite[AM.4]{cartis_adaptive_2011}), we show that DRSOM has a local quadratic convergence and has an \(O(\epsilon^{-3/2})\) complexity to globally converge to a point satisfying the first-order condition \eqref{cond.eps.grad} and the second-order condition in a certain subspace. We identify that this assumption is only needed at the end stage of the algorithm and can be further removed if we periodically perform a \emph{corrector step}, which is a Krylov-like method.
  \item Finally, we perform comprehensive experiments on nonconvex problems. Our results demonstrate its comparable performance to SOMs (e.g. the Newton trust-region method) in terms of the  iteration number and to FOMs (e.g., conjugate gradient method) in terms of the per-iteration cost.  We highlight a set of results in sensor network localization problem, from which we see \drsom{} has the better performance in comparison to a well-known nonlinear conjugate gradient method.
\end{itemize}

Our paper is organized as follows. In \autoref{main.sec.drsom}, we discuss the details of the DRSOM including its important building blocks.  \autoref{main.sec.analysis} is regarding the convergence results of DRSOM such as
the global convergence rate and local convergence rate. Finally, the comprehensive numerical results of DRSOM are presented in \autoref{main.sec.app}.

\section{The DRSOM}\label{main.sec.drsom}
\subsection{Notations and Assumptions}
We first introduce a few notations and technical lemmas. We denote the standard Euclidean norm by $\|\cdot\|$. If $A \in \mathbb{R}^{n \times n}$, then $\|A\|$ is the induced $\mathcal{L}_2$ norm and $A \succeq 0$ indicates that $A$ is positive semidefinite. To facilitate discussion, we denote $g_k=g\left(x_k\right)=\nabla f\left(x_k\right), H_k=$ $H\left(x_k\right)=\nabla^2 f\left(x_k\right)$, and $d_k = x_k - x_{k-1}$ throughout the paper. We use $\mathcal{N}\left(\mu, \sigma^2\right)$ to denote the normal distribution with mean $\mu$, and  variance $\sigma^2$.

In the rest of the paper, we assume $f$ is bounded below and is twice differentiable with $L$-Lipschitz gradient and $M$-Lipschitz Hessian, which is the standard assumption for second-order methods.
\begin{assm}\label{assm.lipschitz}
  $f$ has \(L\)-Lipschitz continuous gradient and \(M\)-Lipschitz continuous Hessian, i.e., for \(\forall x, y \in \rea^n\),
  \begin{equation}\label{eq.assm.lipschitz}
    \|\nabla f(x) - \nabla f(y)\| \le L \|x-y\| \quad \mbox{and}\quad \|\nabla^2 f(x) - \nabla^2 f(y)\| \le M \|x-y\|.
  \end{equation}
\end{assm}
The results in the following lemma are also standard and implied by the first- and second-order Lipshitz continuity.
\begin{lem}[\citet{nesterov_lectures_2018}]\label{lem.lipschitz}
  If \(f:\rea^n \mapsto \rea\) satisfies the \autoref{assm.lipschitz}, then, for all \(x,y\in \rea^n\),
  \begin{subequations}
    \begin{align}
       & |f(y)-f(x)-\nabla f(x)^T(y-x)                                                    \leq \frac{L}{2}\|y-x\|^{2}  \\
       & \left\|\nabla f(y)-\nabla f(x)-\nabla^{2} f(x)(y-x)\right\|                       \leq \frac{M}{2}\|y-x\|^{2} \\
       & \left|f(y)-f(x)-\nabla f(x)^T(y-x)-\frac{1}{2}(y-x)^T\nabla^{2} f(x)(y-x)\right|  \leq \frac{M}{6}\|y-x\|^{3}
    \end{align}
  \end{subequations}
\end{lem}
We will use the results above repetitively throughout the analysis of the paper.

\subsection{Overview of the Algorithm}\label{subsec:overview}
We now give a brief overview of our method. In each iteration of the Dimension-Reduced Second-Order Method (DRSOM), we update
$$
  x_{k+1}=x_k-\alpha_k^1 g_k+\alpha_k^2 d_k,
$$
and the step size $\alpha_k = (\alpha_k^1, \alpha_k^2)$
is determined by solving the following 2-dimensional quadratic model \(m_k(\alpha)\) :
\begin{equation}\label{eq.trs.2d}
  \begin{aligned}
    \min_{\alpha \in\rea^2} & ~m_k(\alpha) := f(\xk) + (c _k)^{T} \alpha+\frac{1}{2} \alpha^{T} \Qk \alpha                                                                 \\
    ~\st~                   & \|\alpha\|_{G_k}:= \sqrt{\alpha^T G_k \alpha} \le \Delta, \quad \mbox{with}\quad G_k=\left[\begin{array}{cc}
                                                                                                                             \left(g_k\right)^T g_k  & -\left(g_k\right)^T d_k \\
                                                                                                                             -\left(g_k\right)^T d_k & \left(d_k\right)^T d_k
                                                                                                                           \end{array}\right],\;
  \end{aligned}
\end{equation}
where
\begin{equation*}
  \Qk =\begin{bmatrix}
    ( \gk)^{T} \Hk \gk  & -( \dk)^{T} \Hk \gk \\
    -( \dk)^{T} \Hk \gk & ( \dk)^{T} \Hk \dk
  \end{bmatrix} \in \mathcal S^{2},\; c _k =\begin{bmatrix}
    -\left\| \gk\right\|^{2} \\
    ( \gk)^{T} \dk
  \end{bmatrix} \in \rea^{2}.
\end{equation*}
The idea of minizing a low-dimensional quadratic model in the objective was previously studied in \citet{ye_second_2005, yuan_subspace_1995}. The novel ingredient in this paper is to impose a $2 \times 2$ trust-region step to determine the step-sizes, which is necessary for solving nonconvex problems. In particular, the trust-region framework is adopted to ensure the global convergence by
controling the progress made at each step. Similar to the standard trust-region method, after computing a trial step $\dkn$, we introduce the reduction ratio $\rho_k$ for $m_k$ at iterate $x_k$,
$$
  \rho_k:=\frac{f\left(x_k\right)-f\left(x_k+d_{k+1}\right)}{m_k(0)-m_k\left(\alpha_k\right)},
$$
if $\rho_k$ is too small, especially when $\rho_k<0$, it means our quadratic model is somehow inaccurate, and it prompts us to decrease the trust region radius; otherwise $m_k$ is good enough, we increase the radius to allow larger search region or keep the radius unchanged.

In the implementation, we also consider the ``Radius-Free" DRSOM by dropping the ball constraint in \eqref{eq.trs.2d} while imposing a quadratic regularization in the objective:
\begin{equation}\label{eq.trs.2d.rf}
  \min_{\alpha \in \rea^2}  ~m_k(\alpha)  + \mu_k\|\alpha\|_{G_k}^2.
\end{equation}
The regularization parameter $\mu_k$ can be tuned to achieve the effect of expanding or contracting the trust region by increasing or decreasing the value of $\mu_k$.

In the following, we present a conceptual DRSOM in \autoref{alg.concept} by incorporating the two alternatives \eqref{eq.trs.2d} and \eqref{eq.trs.2d.rf} to compute the step size $\alpha$, where the adaptive strategy to update the radius of the trust-region constraint or the coefficient of the quadratic regularization is adopted.
%

\IncMargin{11pt}
\setlength{\algomargin}{1pt}
\begin{algorithm}[H]
  \caption{A conceptual DRSOM algorithm}\label{alg.concept}
  \SetAlgoLined
  \KwData{Given $k_{\max}, \beta_1 < 1 < \beta_2, \zeta_1 < \zeta\le 1;\bar{\Delta}>0, \Delta_{0} \in(0, \bar{\Delta})$, and $\eta \in\left[0, \zeta_1\right)$\;}
  \For(){\(k = 1, ..., k_{\max}\)}{
    Solve \eqref{eq.trs.2d} or \eqref{eq.trs.2d.rf} for $\alpha_k$\;
    Compute $\dkn = - \alpha_k^1 \gk +\alpha_k^2 \dk$ and $\rho_k := \frac{f(\xk)-f(\xk+p_k)}{m_k(0)-m_k(p_k)} $\;
    \eIf{$\rho_k>\eta$}{
      Accept the step and update $x_{k+1} = x_{k} + d_{k+1}$
    }{Adjust $\Delta_k $ in \eqref{eq.trs.2d} or $\mu_k $ in \eqref{eq.trs.2d.rf}}
  }
\end{algorithm}
\medskip
In the following, we provide two alternatives to efficiently construct the low-dimensional quadratic model, where the full Hessian matrix is not needed.
\paragraph{Constructing 2-D Quadratic Model Using HVPs}
To efficiently compute the $2 \times 2$ matrix \(Q_k\) in \eqref{eq.trs.2d}, we make use of the decomposition
\[
  Q_k = \left[ -\gk, \dk \right ]^T [-\Hk \cdot \gk , \Hk \cdot \dk ].
\]
Thus, it remains to compute the two Hessian-vector products (HVP) (see \cite{nocedal_numerical_2006}): \(\Hk\gk\) and \( \Hk\dk\).
We adopt the following two strategies to compute those products without request for the full Hessian $\Hk$:
\begin{enumerate}
  \item Finite difference: $\Hk \cdot v \approx \frac{1}{\epsilon} \left[g(x_k+ \epsilon \cdot v) - \gk\right];~ v \in \{\gk, \dk\}$.
  \item Automatic differentiation (AD): $\Hk\gk = \nabla(\frac{1}{2}\gk^T\gk), \Hk\dk = \nabla(\dk^T\gk)$.
\end{enumerate}

\paragraph{Constructing 2-D Quadratic Model Using Interpolation}
Besides HVPs, DRSOM can further benefit from zero-order minimization techniques \cite{conn_global_2009,conn_introduction_2009} and use interpolation to construct the quadratic model.
In particular, given $\gk, \dk$ at iteration $k$, we write the quadratic model:
\begin{equation}\label{eq.interpolation.2nd}
  \begin{aligned}
    m_k(-\alpha_k^1 \gk + \alpha_k^2\dk) ~ & = c_k^T\alpha_k + Q_k \bullet \frac{1}{2}\alpha_k\alpha_k^T \\
    ~                                      & = c_k^T \alpha_k + q_k^Tp_2(\alpha_k),
  \end{aligned}
\end{equation}
where $c_k = [- \|\gk\|^2, \gk^T\dk]$ is known,
$q_k := [Q_{11}, Q_{12}, Q_{22}]^T$
is the vectorized upper-triangular part of $Q_k$, and
\begin{equation}
  \begin{aligned}
    p_2(\alpha_k) & := \left[\frac{(\alpha_k^1)^2}{2}, \alpha_k^1 \alpha_k^2, \frac{(\alpha_k^2)^2}{2}\right]^T
  \end{aligned}
\end{equation}
is the monomial basis of $\alpha_k \in \rea^2$ with order $2$.
We see that if using quadratic interpolation (with exact gradient), the monomial bases are constructed with respect to the \emph{stepsizes} to solve $q_k$ with $3$ unknown entries, which is independent of the problem dimension. To determine the solution, the task left here is to generate an \emph{interpolation set} containing $\ell \ge 3$ different stepsizes, namely $\beta_1, ..., \beta_\ell$.  After that, we will have a linear system forming from $\ell$ extra function evaluations by setting $\alpha_k$ to $\beta_1, ...,\beta_\ell$ respectively. By Taylor expansion and \eqref{eq.interpolation.2nd}, it follows that:
\begin{equation*}
  f(x_k - \alpha^1_k \gk + \alpha_k^2 \dk) = f(x_k) + c_k^T \alpha_k + q_k^Tp_2(\alpha_k) + O(\|- \alpha^1_k \gk + \alpha_k^2 \dk\|^3)
\end{equation*}
Therefore, we solve the following equations with the monomial matrix on the left-hand side:
\begin{align}
  \begin{bmatrix}
    p_2(\beta_1)^T \\
    ...            \\
    p_2(\beta_\ell)^T
  \end{bmatrix}  \cdot q_k & = \begin{bmatrix}
                                 f(x_k - \beta^1_1 \gk + \beta_1^2 \dk) - f(x_k) - c_k^T\beta_1 \\
                                 \cdots                                                         \\
                                 f(x_k - \beta^1_\ell \gk + \beta_\ell^2 \dk) - f(x_k) - c_k^T\beta_\ell
                               \end{bmatrix}.
\end{align}
For computation, the monomial matrix has to be nonsingular. We adopt a simple strategy and randomly pick $\beta_1, ..., \beta_\ell$ over the unit sphere.
In theory, if the iterate moves in the convex hull of sample points, then the error between the function value and the estimation can be uniformly bounded; besides, one can choose more samples so as to alleviate the noise factors \cite{conn_geometry_2008}. In practice, we find that the interpolation method is highly efficient in practice.

In our experience, the interpolation method is preferred over the methods based on finite difference and HVP in most cases in terms of computational speed.


\subsection{Preliminaries of DRSOM}
In this subsection, we present some preliminary analysis of DRSOM.  We introduce the following Lemma regarding the global solution of \eqref{eq.trs.2d} which is widely known.
\begin{lem}\label{lemma:opt-cond}
  The vector \(\alpha^*\) is the global solution to trust-region subproblem \eqref{eq.trs.2d} if it is feasible and there exists a Lagrange multiplier \(\lambda^* \ge 0\) such that $(\alpha^*, \lambda^*)$ is the solution to the following equations:
  \begin{equation}\label{eq.trs.conditions}
    (Q_k + \lambda G_k) \alpha     + c_k   = 0,
    Q_k +  \lambda G_k                    \succeq 0,
    \lambda (\Delta - \|\alpha\|_{G_k})  = 0.
  \end{equation}
\end{lem}
From the construction of \(Q_k\) and \(G_k\), we have that
\begin{equation}
  Q_k + \lambda G_k  = \begin{bmatrix} -\gk^T\\ \dk^T \end{bmatrix} \left(\Hk + \lambda I\right) \begin{bmatrix} -\gk & \dk \end{bmatrix}.
\end{equation}
Therefore, even if \(Q_k\) is indefinite,
there always exists a sufficiently large \(\lambda\) such that condition \eqref{eq.trs.conditions} holds. \citet{ye_new_1991} showed that an \(\epsilon\)-global primal-dual optimizer \((\alpha^*, \lambda^*)\) satisfying \eqref{eq.trs.conditions} can be found in \(O(\log\log(1/\epsilon))\) time. One may also find the optimal solutions by other standard methods in \citet{conn_trust_2000}.
Due to the fact that we only use two directions, the associated subproblems are easy to solve.

We also introduce the normalized problem to facilitate analysis:
\begin{equation}
  \label{eq.trs.2dquadmodelnorm}
  \begin{aligned}
    \min_{\alpha \in\rea^2}~ & f(\xk) + \alpha^T\vk^T\gk+\frac{1}{2} \alpha^{T} \vk^T\Hk\vk \alpha \\
    \st~                     & \|\alpha\| \le \Delta_k,
  \end{aligned}
\end{equation}
where \(V_k\) is the orthonormal bases for \(\spak:= \spa\{\gk, \dk\}\). It is worth mentioning that most of our analysis still goes through for general subspace $\spak$, although the form of $\spak$ is specified in DRSOM.
It is easy to see \eqref{eq.trs.2dquadmodelnorm} and \eqref{eq.trs.2d} are equivalent under a linear transformation.
With a slightly abuse of notation, we do not differentiate $\alpha_k$ and the dual variable $\lambda_k$ to those in \eqref{eq.trs.2d}.
In the next lemma, we note that although problem \eqref{eq.trs.2d} is a 2-dimensional trust region model, its solution can be
transformed into a solution of the ``full-scale'' trust-region problem, as shown below.

\begin{lem}\label{lemma:opt-cond-normal}
  Let $\alpha_k$ and $\lambda_k$ be the solution and the associated Lagrangian multiplier with the trust region constraint to the normalized problem \eqref{eq.trs.2dquadmodelnorm}. Construct \(d_{k+1} = \vk\alpha_k\), then $d_{k+1}$ is a solution to the full-scale problem
  \begin{equation}\label{Prob:t-full-scale}
    \min_{d\in\rea^n} ~\tilde{m}_k(d) := f(\xk) + \gk^Td + \frac{1}{2} d^T\hkap d, ~\st~\|d\| \le \Delta,
  \end{equation}
  such that
  \begin{equation}\label{cond.trs.subspace}
    \begin{aligned}
      (\hkap + \lambda_k I)d_{k+1}  + \gk = 0, ~\hkap + \lambda_k I \succeq 0, ~\lambda_k(\|d_{k+1} \| - \Delta_k) = 0,
    \end{aligned}
  \end{equation}
  where \(
  \hkap = \vk\vk^TH_k\vk\vk^T
  \).
\end{lem}
\begin{proof}
  According to \eqref{eq.trs.conditions}, we have that
  \begin{equation}\label{eq.trs.normal.cond}
    (\vk^T\Hk\vk + I) \alpha_k     + \vk^T\gk   = 0,
    \vk^T\Hk\vk +  \lambda_k I                    \succeq 0,
    \lambda_k (\Delta - \|\alpha_k\|)  = 0.
  \end{equation}
  Multiplying $\vk$ to the left of both sides of the first equation in \eqref{eq.trs.normal.cond} yields that
  \begin{equation*}
    \vk\vk^T\Hk\vk\vk^T\vk\alpha_k + \lambda_k\vk\alpha_k
    = \vk\vk^T\Hk\vk\alpha_k + \lambda_k\vk\alpha_k = -\vk\vk^T\gk
    = -\gk,
  \end{equation*}
  where the first equality is due to \(\vk^T\vk = I\) and the last equality follows from \(\vk\vk^T\) is the projection matrix of $\spak$ and $\gk \in \spak$.
  As a result, we have that:
  \begin{equation*}
    (\vk\vk^T\Hk \vk\vk^T + \lambda_k I) d_{k+1} + \gk=0
  \end{equation*}
  proving the first equation in \eqref{cond.trs.subspace}.
  Due to the second equation in \eqref{eq.trs.normal.cond}, we have that
  $$\alpha^T \vk^T\Hk\vk \alpha + \lambda_k \alpha^T\alpha \ge 0 ,\; \forall\; \alpha \in \rea^2.$$
  By letting $d = \vk \alpha$,
  it is equivalent to
  $$d^T \Hk d + \lambda_k d^Td \ge 0,\; \forall\; d \in \spak$$
  due to \(V_k\) is the orthonormal bases for \(\spak\) and $d^Td=\alpha^T\vk^T \vk \alpha = \alpha^T \alpha $.
  The inequality above is further equivalent to $$\tilde{d}^T(\vk\vk^T\Hk \vk\vk^T + \lambda_k I) \tilde{d} \ge 0,\; \forall \; \tilde{d} \in \rea^n$$
  as \(\vk\vk^T\) is the projection matrix of $\spak$, and this
  proves \(\hkap + \lambda I \succeq 0\) in \eqref{cond.trs.subspace}. Finally, since $d_{k+1}^Td_{k+1}=\alpha_k^T\vk^T \vk \alpha_k = \alpha_k^T \alpha_k $, the last equation in \eqref{cond.trs.subspace} follows from the last equation in \eqref{eq.trs.normal.cond}. Therefore, $d_{k+1}$ is a solution to the problem \eqref{Prob:t-full-scale}.
\end{proof}

The merit of problem \eqref{Prob:t-full-scale} is that we do not require the constraint \(d \in \spak\) and shift the subspace restriction to the matrix \(\hkap\), which can be viewed as the projection of the original Hessian \(\Hk\) in the subspace $\spak$. Therefore, $\hkap$ plays the same role as the approximated Hessian matrix in the inexact Newton method, and DRSOM may similarly be regarded as a cheap quasi-Newton method. In the computations, we actually solve the $2$-dimensional quadratic problem \eqref{eq.trs.2d} for the sake of convenience. However, the facts established for its full-scale counterpart \eqref{Prob:t-full-scale} will be frequently used in the convergence analysis of DRSOM.

\section{Convergence Results}\label{main.sec.analysis}
In this section, we provide a a suite of convergence results of DRSOM. Our analysis includes the equivalence between DRSOM and conjugate gradient (CG) method for strongly convex quadratic programming, the global convergence rate and the local convergence rate of DRSOM for nonconvex setting. Although the subproblem of DRSOM
is constrained to a two dimensional (gradient $g_k$  and momentum $d_k$) subspace,
the global and local convergence rate analysis could be extended to the general
subspace that is not restricted to be
spanned by $g_k$ and $d_k$. Our analysis shows that if the subspace produces a sufficiently good approximated Hessian, DRSOM can achieve extraordinary convergence rate both globally and locally.

\subsection{Finite Convergence for Strictly Convex Quadratic Programming}

We first show DRSOM has finite convergence for strongly convex quadratic programming.
\begin{equation}\label{eq.conv.qp}
  \min f(x) = \frac{1}{2}x^TAx + a^Tx,
\end{equation}
where $A \succ 0$. In this case, we drop the the trust-region constraint for DRSOM, i.e., \(\Delta_k\) is sufficiently large, \(\lk = 0\) for all \(k\). We have the following theorem.
\begin{thm} \label{thm.conv.qp.equiv}
  If we apply DRSOM to \eqref{eq.conv.qp} with no radius restriction, i.e., \(\Delta\) is sufficiently large, then the DRSOM generates the same iterates of conjugate gradient method, if they start at the same point \(x_0\).
\end{thm}
\begin{proof}
  To show its equivalence to the conjugate gradient method, we only have to prove the iterate \(\xk\) by DRSOM minimizes \(f(x) \) in the subspace such that:
  \[\xk \in \spak =  x_0 + \spa\{d_1, ..., d_k\}.\]
  Since there is no radius, the solution that minimizes \(m_k\) strictly corresponds to the optimizer of \(f(x)\) in the subspace \(\xk + \spa\{\gk, \dk\}\). In other words, the iterate of DRSOM can also be retrieved by simply choosing the stepsizes such that \(\xkn = \arg\min f(x)\) and \(\xkn = \xk - \alpha_k^1 \gk + \alpha_k^2 \dk\).
  From such a perspective, we show that the \(\xk\) is equivalent to \(\tilde x_k\) by the conjugate gradient method.

  Note \(\tilde \xk \) minimizes \(f(x)\) over the subspace below:
  \begin{equation*}
    x_0 + \spa\{\tilde d_1, ..., \tilde \dk\},
  \end{equation*}
  where \(\tilde d_1, ..., \tilde \dk\) are conjugate directions for CG. By construction, we see \(x_1 = x_0 + \alpha^1_0 g_0\) and \(d_1 = \alpha^1_0 g_0\), so that \(x_1 = \tilde x_1\).
  Assume it holds for \(k\), we know for CG:
  \[\tilde  \gk \in  \spa\{\tilde \dk, \tilde d_{k+1}\},\]
  and since \(\gk = \tilde \gk, \dk = \tilde \dk\), we have
  \begin{equation}
    \spa\{\tilde \dk, \tilde d_{k+1}\} = \spa\{\dk, \gk\},
  \end{equation}
  Now we know from the next minimizer \(\xkn \in \xk + \spa\{\dk, \gk\}\):
  \[\xkn \in \tilde x_k + \spa\{\tilde \dk, \tilde d_{k+1}\},\]
  and since \(\tilde x_{k+1}\) minimizes \(f\) over both \(x_0 + \spa\{d_1, ..., \dk, \tilde d_{k+1}\}\) and \(\xk + \spa\{\tilde d_{k+1}\}\), we have that \(\xkn = \tilde x_{k+1}\) as the desired result.
\end{proof}
The equivalence of quadratic minimization over $\spak$ and conjugate gradient method was also established in \citet{yuan_subspace_1995}. We provide the proof of Theorem \ref{thm.conv.qp.equiv} for the completeness of the paper.

\subsection{Global Convergence Rate}

To analyze the convergence rate of DRSOM, we need to make the following assumption regarding the approximated Hessian $\hkap$, which is commonly used in the literature; see \cite{dennis_quasi-newton_1977,cartis_adaptive_2011-1,curtis_trust_2017,xu_newton-type_2020}.
\begin{assm}\label{assm.hessian}
  The approximated Hessian matrix \(\hkap\) along subspace \(\mathcal L_k\) satisfies:
  \begin{equation}\label{eq.assm.hessian}
    \|(\Hk - \hkap)\dkn\| \le  C \|\dkn\|^2
  \end{equation}
\end{assm}
Although we adopt the adaptive strategy to choose the radius in the implementation, our analysis is conducted under the fixed-radius strategy such that a step is always accepted for simplicity.
In terms of the global convergence rate, we show that DRSOM has an \(O\left(\epsilon^{-3/2}\right)\) complexity to converge to the first-order stationary point and second-order point approximately in the subspace.
Let us assume that \autoref{assm.hessian} holds for the following analysis.

To prove our main theorem, we need to establish the following lemma as preparation.
\begin{lem}[Model reduction]\label{lemma:mk-reduc}
  At iteration \(k\),
  let $d_{k+1}$ and $\lambda_k$ be the solution and Lagrangian multiplier
  constructed in \autoref{lemma:opt-cond-normal}.
  If \(\lk > 0\) , we have the following amount of decrease on $\tilde{m}_k$:
  \begin{equation}
    \tilde{m}_k(d_{k+1}) - \tilde{m}_k(0) = - \frac{1}{2} \lambda_k \Delta_k^2.
  \end{equation}
\end{lem}

\begin{proof}
  In view of \autoref{lemma:opt-cond-normal}, the optimal condition \eqref{cond.trs.subspace} holds.
  Then, we have \(\|\dkn\| = \Delta_k\) due to \(\lk > 0\) and  that:
  \begin{equation}
    \begin{aligned}
      \tilde{m}_k(\dkn) - \tilde{m}_k(0) & = \gk^T\dkn + \frac{1}{2} \dkn^T\hkap \dkn                                 \\
                                         & = -\frac{1}{2} \dkn^T  (\hkap + \lk I) \dkn + \frac{1}{2} \dkn^T\hkap \dkn \\
                                         & = - \frac{1}{2} \lk \Delta_k^2,
    \end{aligned}
  \end{equation}
  which completes the proof.
\end{proof}

Taking \(\Delta_k = \Delta =  \frac{2\sqrt{\epsilon}}{M}\) and combining with the reduction of the quadratic model, we conclude that DRSOM generates sufficient decrease at every iteration \(k\) as long as \(\lambda_k \geq \sqrt{\epsilon}\). The following analysis based on a \emph{fixed} trust-region radius is motivated from \citet{luenberger_linear_2021}.

\begin{lem}[Sufficient decrease]\label{lem: function value decrease}
  At iteration \(k\), take \(\Delta_k = \Delta =  \frac{2\sqrt{\epsilon}}{M}\), and let $d_{k+1}$ and $\lambda_k$ be the solution and Lagrangian multiplier
  obtained in \autoref{lemma:opt-cond-normal}.
  If \(\lk \ge \sqrt{\epsilon}\) , we have the following amount of function value decrease,

  \begin{equation}\label{eq.decrease}
    f(\xkn) \le f(\xk) - \frac{2}{3M^2} \epsilon^{3/2}.
  \end{equation}
\end{lem}

\begin{proof}
  Since $d_{k+1} \in \spak $ and \(\vk\vk^T\) is the projection matrix of $\spak$, it holds that $$d_{k+1}^T \hkap d_{k+1}=d_{k+1}^T \vk\vk^T\Hk \vk\vk^T d_{k+1} = d_{k+1}^T \Hk  d_{k+1}.$$
  Moreover,
  with second-order Lipschitz continuity and tailor expansion, we immediately have:
  \begin{equation}\label{eq.reduction}
    \begin{aligned}
      f(\xkn) & \le f(\xk) + (\gk)^T\dkn + \frac{1}{2}(\dkn)^T\Hk(\dkn) + \frac{M}{6}\|\dkn\|^3       \\
              & = f(\xk) + (\gk)^T\dkn + \frac{1}{2}(\dkn)^T\tilde{H}_k(\dkn) + \frac{M}{6}\|\dkn\|^3 \\
              & =  f(\xk) - \frac{1}{2} \lambda_k \Delta^2 + \frac{1}{6} M \Delta^3                   \\
              & = f(\xk) - \frac{2\lambda_k \epsilon}{M^2} + \frac{4\epsilon^{3/2}}{3M^2}
    \end{aligned}
  \end{equation}
  where the second last equality is due to \autoref{lemma:mk-reduc} and \(\|\dkn\| = \Delta_k = \Delta\) from the optimality condition
  \eqref{cond.trs.subspace}
  with \(\lk > 0\). Noting that \(\lambda_k \ge \sqrt \epsilon\), inequality \eqref{eq.reduction} implies that
  \begin{align*}
    f(\xkn) & \le f(\xk) - \frac{2}{3M^2} \epsilon^{3/2}.
  \end{align*}
\end{proof}
The following result states that when \(\lk \le \sqrt{\epsilon}\) and the Hessian regularity condition holds, we can terminate the process at the next iterate \(\xkn\) that approximated satisfies the first-order condition and the second-order condition in the subspace.
\begin{lem}\label{lemma:global-small-lambda}
  At iteration \(k\), if the Lagrangian multiplier $\lk$ associated with the trust region constraint in \eqref{eq.trs.2d} satisfies \(\lk \leq \sqrt{\epsilon}\)  and Hessian regularity condition \eqref{eq.assm.hessian} holds, then the iterate \(\xkn\) approximately satisfies the first-order
  condition, and the second-order condition in the subspace \(\mathcal{L}_k\).
\end{lem}
\begin{proof}
  Suppose $d_{k+1}$ is the solution
  obtained in \autoref{lemma:opt-cond-normal}. By second-order Lipschitz continuity and the first equation in the optimality condition \eqref{cond.trs.subspace}, we have that:
  \begin{equation}
    \begin{aligned}
      \|\gkn\| & \le \|\gkn - \gk - \Hk \dkn\| + \|(\Hk - \hkap)\dkn\| + \|(\gk + \hkap \dkn)\|
      \\
               & \le \left\|\int_0^1 \left[ \nabla^2 f(\xk + \tau \dkn) - \Hk\right] \dkn d\tau \right\| + \|(\Hk - \hkap)\dkn\|+ \lk \|\dkn\| \\
               & \le \frac{1}{2} M \|\dkn\|^2 + \lk \|\dkn\| + \|(\Hk - \hkap)\dkn\|
    \end{aligned}
  \end{equation}
  In view of  \autoref{assm.hessian}, \(\lk \leq \sqrt{\epsilon}\) and $\|\dkn\| \le \Delta =  \frac{2\sqrt{\epsilon}}{M}$, we immediately have
  \begin{equation}\label{eq.global.grad}
    \begin{aligned}
      \|\gkn\| & \le \left(\frac{1}{2} M + C\right) \|\dkn\|^2 + \lk \|\dkn\|                   \\
               & \le \left( \frac{1}{2}M + C\right) \frac{4\epsilon}{M^2} + \frac{2\epsilon}{M} \\
               & \le \left(\frac{4}{M} + \frac{4C}{M^2}\right) \epsilon.
    \end{aligned}
  \end{equation}
  As for the second-order condition, the second condition in \eqref{cond.trs.subspace} and $\lk \leq \sqrt{\epsilon}$ imply that
  \begin{align}
    \nonumber -\sqrt{\epsilon}I\preceq -\lambda_k I\preceq \hkap & = V_k V_k^T H_{k+1} V_k V_k^T + \hkap - V_k V_k^T H_{k+1} V_k V_k^T            \\
    \nonumber                                                    & = V_k V_k^T H_{k+1} V_k V_k^T + \vk\vk^T (\Hk-H_{k+1}) \vk\vk^T                \\
    \nonumber                                                    & \preceq V_k V_k^T H_{k+1} V_k V_k^T + \|\vk\vk^T (\Hk-H_{k+1}) \vk\vk^T\| I    \\
    \nonumber                                                    & \preceq V_k V_k^T H_{k+1} V_k V_k^T + \|\vk\vk^T\|\|H_{k+1}-\Hk\|\|\vk\vk^T\|I \\
    \nonumber                                                    & = V_k V_k^T H_{k+1} V_k V_k^T + \| H_{k+1}-\Hk\|I                              \\
    \nonumber                                                    & \preceq V_k V_k^T H_{k+1} V_k V_k^T + M \|d_{k+1}\|I                           \\
                                                                 & \preceq V_k V_k^T H_{k+1} V_k V_k^T + 2\sqrt{\epsilon} I,
  \end{align}
  where the second last matrix inequality is due to the Lipschitz continuity of the Hessian and the last matrix inequality follows from $\|\dkn\| \le \Delta =  \frac{2\sqrt{\epsilon}}{M}$. Hence, it holds that
  \begin{equation*}
    V_k V_k^T H_{k+1} V_k V_k^T  \succeq -3\sqrt{\epsilon} I,
  \end{equation*}
  which indicates that  \(H_{k+1}\) is approximately positive semi-definite  in the subspace \(\spak\).
\end{proof}
\begin{thm}\label{thm.global}
  Under the fixed-radius strategy by setting \(\Delta_k =  \frac{2\sqrt{\epsilon}}{M}\),
  the DRSOM runs at most \(O\left(\frac{3}{2}M^2(f(x_0) - f_{\inf})\epsilon^{-3/2}\right)\) iterations to reach an iterate \(\xkn\) that satisfies the first-order condition \eqref{cond.eps.grad}, and
  the approximated second-order condition in the subspace \(\spak\):
  \(
  \vk\vk^TH_{k+1}\vk\vk^T    \succeq -3\sqrt{\epsilon}I,
  \)
  where $\vk$ is the orthonormal bases for \(\spak\).
\end{thm}
\begin{proof}
  According to \autoref{lemma:global-small-lambda}, when $\lk \le \sqrt{\epsilon}$, we already obtain
  an iterate that approximately satisfies the first-order
  condition, and the second-order condition in certain subspace.
  On the other hand, when \(\lk > \sqrt{\epsilon}\),
  Lemma \ref{lem: function value decrease} indicates that the objective function has a amount of decrease \(\frac{2}{3M^2} \epsilon^{3/2}\) at every iteration \(k\). Note that the total amount of decrease cannot exceed \(f(x_0) - f_{\inf}\). Therefore, the number of iterations with \(\lk > \sqrt{\epsilon}\) is upper bounded by
  \begin{equation*}
    O\left(\frac{3}{2}M^2(f(x_0) - f_{\inf})\epsilon^{-3/2}\right),
  \end{equation*}
  which thus is also the iteration bound of our algorithm.
\end{proof}
\subsection{Local Convergence Rate}
Regarding the local rate of convergence, we can prove that DRSOM have a local quadratic convergence rate. To this end, we first present some technical lemmas.
\begin{lem}\label{lem.subspace.drsom-newton}
  Let \(V_k\) be the orthonormal bases for \(\spak\), suppose \(\ak\) is the solution to the normalized problem \eqref{eq.trs.2dquadmodelnorm}.
  Let \(\dkn = \vk\alpha_k\), then the following inequality holds:
  \begin{equation}\label{prob:subspace-Newton}
    \|\dkn^{SN} - \dkn\| \le \frac{1}{\mu }\lk \|\dkn\|, \\
  \end{equation}
  where $\lambda_k$ is the Lagrangian multiplier associated with the trust region constraint in  \eqref{eq.trs.2dquadmodelnorm}, and $d_{k+1}^{SN} = V_k \alpha_{k}^{SN}$ is the subspace Newton step with $\alpha_{k}^{SN}$ defined by:
  \begin{equation}\label{eqn:alpha-SN}
    \alpha_k^{SN} = \arg\min_{\alpha\in\rea^2}~ f(\xk) + \alpha^T\vk^T\gk+\frac{1}{2} \alpha^{T} \vk^T\Hk\vk \alpha.
  \end{equation}
\end{lem}

\begin{proof}
  Since \(\ak\) is a solution to the normalized problem \eqref{eq.trs.2dquadmodelnorm}, the optimality condition gives that
  $$
    (\vk^T\Hk\vk + \lk I)\alpha_k = -\gk^T\vk.
  $$
  As \(\ak^{SN}\) is a solution to problem \eqref{prob:subspace-Newton}, due to optimality condition it holds that
  $$
    \vk^T\Hk\vk \alpha_k^{SN} = -\gk^T\vk.
  $$
  Combining the above two equations yields that
  \begin{equation}\label{eqn:ak-akSN}
    \vk^T\Hk\vk (\alpha_k - \ak^{SN})  = -\lk \alpha_k.
  \end{equation}
  Moreover, note that
  \begin{equation}\label{eqn:norm-invar}
    \mbox{for any}\; \alpha \neq 0, \;\mbox{we have}\; \vk \alpha \neq 0 \;\mbox{and}\; \| \vk \alpha \| = \|\alpha\|,
  \end{equation}
  which combined with $\Hk \succeq \mu I$ implies that
  \begin{equation*}
    \alpha ^T \vk^T \Hk \vk \alpha \geq \mu \| \vk \alpha\|^2 = \mu \| \alpha \|^2.
  \end{equation*}
  Therefore $\vk^T \Hk\vk \succeq \mu I_2$ holds (also implies that $\vk^T \Hk\vk$ is nonsingular), it follows that
  \begin{equation}
    \| \ak^{SN} - \ak\| \leq\|(\vk^T\Hk\vk)^{-1}\| \|\vk^T\Hk\vk (\ak^{SN} - \ak)\| \leq \frac{1}{\mu}\lk \|\ak\|,
  \end{equation}
  where the second inequlity is due to \eqref{eqn:ak-akSN}.
  Therefore, by combining the construction of $\dkn^{SN}$ and $\dkn$ with \eqref{eqn:norm-invar} we conclude that
  \begin{equation}
    \|\dkn^{SN} - \dkn\| = \|\ak^{SN} - \ak\| \le \frac{1}{\mu }\lk \|\ak\| =\frac{1}{\mu }\lk \|\dkn\|.
  \end{equation}
\end{proof}
We now ready to provide the following key result to analyze the local convergence rate of our algorithm, where we assume
it converges to a strict local optimum \(x^*\) such that \(H(x^*) \succeq \mu I\) for some $\mu > 0$.
\begin{lem}
  Suppose the iterate of DRSOM $\xk$ converges to $x^*$ which satisfies \(H(x^*) \succeq \mu I\), when $\xk$ is sufficiently close to $x^*$, then we have:
  \begin{equation}\label{eq.local.superlinear}
    \| \xkn-x^*\|\leq
    \frac{M}{\mu}\|\xk-x^*\|^2
    + \frac{1}{\mu}\|(\Hk-\hkap)\dkn\|
    + \left(\frac{2L}{\mu^2} + \frac{1}{\mu}\right) \lk \|\dkn\|.
  \end{equation}
\end{lem}
\begin{proof}
  We first write
  \begin{equation}\label{eqn:xk-xstar}
    \begin{aligned}
      \|{x_{k+1}-x^*}\| & = \|\xk +\dkn-x^*\|                                                                    \\
                        & \leq \| \xk+d_{k+1}^N-x^*\|+ \| d_{k+1}^{SN} - d_{k+1}^N\|+\| d_{k+1}- d_{k+1}^{SN}\|,
    \end{aligned}
  \end{equation}
  where $d_{k+1}^N = -H_k^{-1}g_k$ is the standard Newton step and $d_{k+1}^{SN} = V_k \alpha_{k}^{SN}$ is the subspace Newton step with $\alpha_{k}^{SN}$ defined by
  \eqref{eqn:alpha-SN}.
  The first term in \eqref{eqn:xk-xstar} is upper bounded by \(\frac{M}{\mu}\|\xk-x^*\|^2\) due to the standard analysis in Newton's method.
  To bound the second term, we note that $\alpha_k^{SN} = (V_k^T H_k V_k)^{-1}V_k^Tg_k$ with $V_k^TV_k = I$ as it is a solution to problem \eqref{eqn:alpha-SN}, which implies that
  $$
    \begin{aligned}
      \tilde{H}_k d_{k+1}^{SN} & = \tilde{H}_k V_k (V_k^T H_k V_k)^{-1}V_k^Tg_k               \\
                               & = V_k V_k^T H_k V_k V_k^{T} V_k (V_k^T H_k V_k)^{-1}V_k^Tg_k \\
                               & = V_k V_k^T g_k = g_k,
    \end{aligned}
  $$
  where the last equality is due to $V_k V_k^T$ is the projection matrix of the subspace $\spak$ and $g_k \in \spak$.
  Then, the second term can be further bounded above as follows
  \begin{align}
    \nonumber     \| d_{k+1}^{SN}- d_{k+1}^N \| & = \| d_{k+1}^{SN}+ H_k^{-1}g_k\|                                                                                \\
    \nonumber                                   & = \| d_{k+1}^{SN}-H_k^{-1}\tilde{H}_kd_{k+1}^{SN}\|                                                             \\
    \nonumber                                   & = \| H_k^{-1} (H_k-\tilde{H}_k)d_{k+1}^{SN}\|                                                                   \\
    \nonumber                                   & \leq \| H_k^{-1}\| \| (H_k-\tilde{H}_k)d_{k+1}^{SN}\|                                                           \\
    \nonumber                                   & \leq \frac{1}{\mu}\| (H_k-\tilde{H}_k)d_{k+1}^{SN}\|                                                            \\
    \nonumber                                   & \leq \frac{1}{\mu} \left(  \| (H_k-\tilde{H}_k)d_{k+1}\|+ \| (H_k-\tilde{H}_k)(d_{k+1}^{SN}-d_{k+1})\| \right)  \\
                                                & \leq \frac{1}{\mu} \|  (H_k-\tilde{H}_k)d_{k+1}\|+\frac{1}{\mu}\| H_k-\tilde{H}_k\| \|(d_{k+1}^{SN}-d_{k+1})\|,
  \end{align}
  Combining the above inequalities, we have that
  \begin{equation}
    \begin{aligned}
          & \|{x_{k+1}-x^*}\|                                                                                                                 \\
      \le & \frac{M}{\mu}\|\xk-x^*\|^2 + \frac{1}{\mu}\|(\Hk-\hkap)\dkn\| + \left(\frac{1}{\mu} \|\Hk-\hkap\| + 1\right) \|\dkn - \dkn^{SN}\| \\
      \leq
          & \frac{M}{\mu}\|\xk-x^*\|^2
      + \frac{1}{\mu}\|(\Hk-\hkap)\dkn\|
      + \left(\frac{2L}{\mu^2} + \frac{1}{\mu}\right) \lk \|\dkn\|,
    \end{aligned}
  \end{equation}
  where the last inequality follows from \autoref{lem.subspace.drsom-newton} and $\max\{ \|H_k\|, \|\tilde{H}_k\| \} \le L$ due to the Lipschitz continuity of the gradient.
\end{proof}
\begin{thm}\label{thm.local}  Suppose \(x^*\) is a second-order stationary point such that \(H(x^*) \succeq \mu I\) for some $\mu > 0$.
  Then if $\xk$ is sufficiently close to $x^*$, DRSOM converges to \(x^*\) quadratically, namely:
  $
    \|{x_{k+1}-x^*} \| \le O(\|\xk - x^*\|^{2}).
  $
\end{thm}
\begin{proof}
  It suffices to further upper bound \eqref{eq.local.superlinear}.
  We only consider the scenario that
  \(\lk \le \sqrt{\epsilon}\), as otherwise, the objective function has an amount of decrease \(\frac{2}{3M^2} \epsilon^{3/2}\) at every iteration \(k\) by Lemma \ref{lem: function value decrease} and this will occur in a very limited time when $x_k$ is sufficiently close to $x^*$.

  Since \eqref{eq.assm.hessian} holds, one has that
  $$
    \frac{1}{\mu}\|(\Hk-\hkap)\dkn\| \le \frac{C}{\mu} \|\dkn\|^2.
  $$
  Moreover, as we adopt the fixed radius strategy, $\lambda_k = 0 $
  whenever $\|d_{k+1}\| < \Delta$. In the case of $0<\lk \leq \sqrt{\epsilon}$, we have \( \lk \leq \sqrt{\epsilon} = \frac{M}{2} \Delta = \frac{M}{2}\|\dkn\|\). Therefore, in both cases, we have $\lambda_k \|d_{k+1}\| \le \frac{M}{2} \|\dkn\|^2$. Consequently, \eqref{eq.local.superlinear} can be bounded above by
  \begin{equation*}
    \|{x_{k+1}-x^*}\| \le \frac{M}{\mu}\|\xk-x^*\|^2 + \frac{C}{\mu} \|\dkn\|^2 + \left(\frac{2L}{\mu^2} + \frac{1}{\mu}\right) \frac{M}{2} \|\dkn\|^2.
  \end{equation*}
  Note that
  \begin{equation}
    \begin{aligned}
      \|\dkn\| & \leq \|{x_{k}-x^*+\dkn}\| + \|x_k - x^*\|                                                                                                         \\
               & = \|{x_{k+1}-x^*}\| + \|x_k - x^*\|                                                                                                               \\
               & \leq \frac{M}{\mu}\|\xk-x^*\|^2 + \frac{C}{\mu} \|\dkn\|^2 + \left(\frac{2L}{\mu^2} + \frac{1}{\mu}\right) \frac{M}{2} \|\dkn\|^2 + \|x_k - x^*\| \\
               & \leq \frac{M}{\mu}\|\xk-x^*\|^2 + \|x_k - x^*\| + O(\|\dkn\|^2). \notag
    \end{aligned}
  \end{equation}
  By rearranging the terms, we have
  \begin{equation}\label{eq.dkn}
    \| \dkn \|-O(\|\dkn\|^2)\leq \frac{M}{\mu}\| \xk-x^*\|^2+\|\xk-x^*\|.
  \end{equation}
  From the assumption $\xk$ converges to $x^*$, it holds that $\dkn\rightarrow 0$, when $k$ is sufficiently large. Thus, inequality \eqref{eq.dkn} implies $\|\dkn\|\leq \|\xk-x^*\|$, i.e. \(\|\dkn\| = O(\|\xk - x^*\|)\), which in return shows that \(\|{x_{k+1}-x^*}\| \le O(\|\xk - x^*\|^2)\).
\end{proof}

\subsection{A Corrector Step}

In fact, the inequality \eqref{eq.assm.hessian} in \autoref{assm.hessian} plays a crucial role in our convergence analysis. To ensure \eqref{eq.assm.hessian}, a popular strategy is to apply the Lanczos method; see \cite{cartis_adaptive_2011-1,curtis_trust_2017,conn_trust_2000}. Although we can not rigorously establish the validness of \autoref{assm.hessian} using merely $\gk$ and $\dk$, we manage to find out that it is not always required in the iteration process. In view of \eqref{eq.reduction}, we see that when $\lambda_k > \sqrt{\epsilon}$, the iterate has sufficient decrease in function value and thus \eqref{eq.assm.hessian} does not really matter. Namely, the objective function has a amount of decrease \(\frac{2}{3M^2} \epsilon^{3/2}\) (by Lemma \ref{lem: function value decrease}). It is only when $\lambda_k \le \sqrt{\epsilon}$ that DRSOM may reach a stationary point, then one must check for the sake of \eqref{eq.global.grad}. This fact also applies to the last part of the proof for \autoref{thm.local}.

From these observations, we restrict our attention to the case where \(\lk \le \sqrt{\epsilon}\) and \eqref{eq.assm.hessian} does not hold. Before we delve into a detailed discussion of our treatment, we first provide two different viewpoints of the assumption itself. Firstly, the left-hand side of \autoref{assm.hessian} can be interpreted as the \textit{approximated} residual of Newton's equations using an inexact Newton method \cite{dembo_inexact_1982}. In general, for the Newton step, the equation $H_k d_{k+1}^{N} = -g_k$ is valid. However, for the inexact Newton step $\bar d_{k+1}$, the equation no longer holds and hence we define the corresponding residual as $r_k = \|H_k \bar d_{k+1} + g_k\|$. The next Lemma demonstrates that the left-hand side of \autoref{assm.hessian} is very close to the inexact newton step residual. Since this fact comes from very early literature, we include it here for completeness.
\begin{lemma}
  Define the inexact Newton step residual $r_k:=\|H_k d_{k+1} + g_k\| \in \rea$ , then the LHS of \autoref{assm.hessian} has the following relation with the residual $r_k$ satisfying
  \begin{equation}
    \label{ineq.relation rk lhs}
    r_k - \|(\Hk -\tilde{H}_k)\dkn\| \leq O(\epsilon).
  \end{equation}
  Moreover, if the trust-region constraint is inactive, we have
  \begin{equation}
    \label{eq.relation rk lhs 2}
    r_k = \|(\Hk -\tilde{H}_k)\dkn\|.
  \end{equation}
\end{lemma}
\begin{proof}
  Recalling the results in \autoref{lemma:opt-cond-normal}, $d_{k+1}$ satisfies
  \begin{equation*}
    \tilde{H}_k d_{k+1} + \lambda_k d_{k+1} = -g_k.
  \end{equation*}
  Substituting the above equation into the definition of $r_k$, we obtain
  \begin{equation}
    \begin{aligned}
      r_k & =    \|H_k d_{k+1} + g_k\|                                     \\
          & =    \|H_k d_{k+1} - \tilde{H}_k d_{k+1} - \lambda_k d_{k+1}\| \\
          & \leq \|(\Hk -\tilde{H}_k)\dkn\| + \lambda_k \|d_{k+1}\|.
    \end{aligned}
  \end{equation}
  Note that the \autoref{assm.hessian} is only required when $\lambda_k \leq \sqrt{\epsilon}$, and $\|d_{k+1}\| \leq \Delta_k = 2\sqrt{\epsilon}/M$, therefore the inequality \eqref{ineq.relation rk lhs} holds. As for the second statement, if the trust-region constraint is inactive, then we have $\lambda_k=0$, and \eqref{eq.relation rk lhs 2} follows immediately.
\end{proof}
The seminal work \cite{dembo_inexact_1982} addresses the conditions and properties of solving a Newton equation \emph{inexactly}.  More recently, \citet{gould_error_2020} take residual as $\|\Hk\dkn + \gk + \lambda_k \dkn\|$, which can be applied to trust-region and $p$-th order regularized subproblems. These facts and analyses justify the iterative methods with evolving subspaces.

Our second viewpoint is that the left-hand side of \autoref{assm.hessian} is equivalent to a maximization problem, whose solution corresponds to the ``worst'' direction. That is, our current projected Hessian has the largest error along such a direction, which is shown below.

\begin{lemma}\label{lem.corrector.interpretation}
  Suppose $\dkn$ is calculated over a subspace $\mathcal L_k$, then we have for the LHS of \eqref{eq.assm.hessian},
  \begin{equation}\label{compute u}
    \|(H_k - \tilde \Hk)d_{k+1}\| = \underset{\|u\|= 1, u \perp \mathcal L_k}{\max} \quad u^T\Hk d_{k+1}.
  \end{equation}
\end{lemma}
\begin{proof}
  Let $\vk$ be the orthonormal basis for $\mathcal L_k$, then,
  \begin{subequations}
    \begin{align}
      \|(H_k - \tilde \Hk)d_{k+1}\| = & ~ \|(\Hk - V_k V_k^T\Hk \vk V_k^T)d_{k+1}\|                               \\
      \label{eq.sse.proj}=            & ~ \|(\Hk - V_k V_k^T\Hk) d_{k+1}\|                                        \\
      =                               & ~ \| (I - \vk\vk^T) \Hk \dkn \|                                           \\
      =                               & ~  \underset{\|x\|= 1}{\max} ~ x^T(I - \vk\vk^T) \Hk d_{k+1}              \\
      \label{eq.sse.moti}=            & ~      \underset{\|u\|= 1, u \perp \mathcal L_k}{\max} ~ u^T \Hk d_{k+1}.
    \end{align}
  \end{subequations}
  Similar to previous analyses, \eqref{eq.sse.proj} follows from $\dkn \in \mathcal L_k$.  The last line \eqref{eq.sse.moti} is because that $(I - \vk\vk^T)$ is also an orthonormal projection matrix.
\end{proof}
In the ideal case of $\Hk\dkn \in \mathcal L_k$, \eqref{eq.sse.moti} equals to 0; otherwise, $\Hk d_{k+1}$ is the the solution to \eqref{eq.sse.moti} and can be viewed as the \emph{steepest} direction.
Furthermore, if we hope to minimize the optimal value of \eqref{eq.sse.moti} by
enlarging the subspace $\mathcal L_k$ with a new direction, a Krylov direction $\Hk\dkn$ based on the current iterate is the best choice. This motivates us to consider the following subproblem with expanded subspace:
\begin{equation}\label{eq.subj}
  \begin{aligned}
    \min_{d\in\mathbb R^n}~ & m_k(d)                                                                       \\
    ~\st{}~                 & d\in \mathcal L_k^{j}:=\text{span}\{\mathcal L_k^{j-1}, H_k d_{k+1}^{j-1}\}.
  \end{aligned}
\end{equation}
In this view, we introduce a corrector step in \autoref{alg.corrector} to iteratively introduce Krylov-type directions and solve the resulting subproblem \eqref{eq.subj} in a larger subspace until \autoref{assm.hessian} is met.

\IncMargin{11pt}
\setlength{\algomargin}{1pt}
\begin{algorithm}[H]
  \small
  \caption{The Corrector Step}\label{alg.corrector}
  \SetAlgoLined
  \KwData{Given iterate $k$, $\gk, \Hk, \dk$, $d_{k+1}^0 = \dkn$\;}
  \For(){\(j = 1, ..., n\)}{
    Solve new problem \eqref{eq.subj} and let $d^{j}_{k+1} , \lambda^{j}_{k}$ be primal-dual solution\;
    \If{$\lambda_{k+1}>\sqrt{\epsilon}$ or \eqref{eq.assm.hessian} holds}{
      Return $\dkn := d^{j}_{k+1},~\lambda_k:=\lambda^{j}_{k}$}
  }
\end{algorithm}

Although \eqref{eq.subj} utilizes a new Krylov direction at each iteration, we should emphasize that it is not a Krylov subspace method. Indeed, a \emph{genuine} Krylov subspace method  of order $j$ constructs the following subspace,
$$\mathcal K^j(\Hk, \gk):=\operatorname{span}\left\{\gk, \Hk \gk, \ldots, \Hk^{j-1} \gk\right\},$$
where the vector $\gk$ always appears in the matrix-vector products, while
the matrix-vector products used in \eqref{eq.subj} of \autoref{alg.corrector} do not use a constant vector and it is taken as $\dkn^0, \dkn^1,...$ and so forth. The crux of our approach is that we directly tackle LHS of \eqref{eq.assm.hessian} that may gain an advantage over traditional Krylov methods.
To justify \autoref{alg.corrector}, we provide the following results.

\begin{lem}[Finite convergence of \autoref{alg.corrector}]
  If the corrector step starts with $\mathcal L_k^0 = \mathrm{span}\{\gk, \dk\}$, then \autoref{alg.corrector} takes at most $n-1$ steps to stop.
\end{lem}

\begin{proof}
  The proof immediate follows
  if $\Hk d^j_{k+1} \in \mathcal L_k^j$, as in this case  $u^T \Hk d_{k+1} = 0$ holds for any feasible $u$ in \eqref{eq.subj}, and \eqref{eq.assm.hessian} must hold. Moreover, the optimal value of \eqref{eq.subj} also equals to $0$ if $\mathcal L_k^{j-1}$ spans the whole space $\rea^n$. Since the dimension of the subspace  $\mathcal L_k^{j}$ in the next iteration is increased by one if \autoref{alg.corrector} is not stopped, and it starts with a two dimensional subspace, \autoref{alg.corrector}
  terminates within at most $n-1$ steps.
\end{proof}
We would like to remark that
the finite iteration bound of
\autoref{alg.corrector}
is very conservative. In practice, since we only require \eqref{eq.assm.hessian} that is a relaxation of $\|(\Hk - \tilde \Hk) \dkn\| = 0$ exactly, \autoref{alg.corrector} usually consumes much less iterations than the upper bound $n-1$.


However, it is possible that the above method in return produces a larger \(\lambda_{k} > \sqrt{\epsilon}\), in which case we terminate the inner loop at some $j \ge 1$ of the corrector method \autoref{alg.corrector} and proceed the main loop at $\xkn$ in \autoref{alg.concept}. This again aligns with our discussion on \eqref{eq.reduction} above.
As we apply the corrector step \emph{periodically} whenever \(\lambda_k \le \sqrt{\epsilon}\), for every iteration with \(\lambda_k \ge \sqrt{\epsilon}\),  Thus the number of corrector steps is very limited when $x_k$ is close to convergence. Finally, we terminate the algorithm as soon as \(\lk \le \sqrt{\epsilon}\) and  \eqref{eq.assm.hessian} hold simultaneously.



\section{Numerical Experiments}\label{main.sec.app}

In this section, we provide detailed experiments of DRSOM on nonconvex optimization problems. 
We first run DRSOM on the \(\mathcal L_2 - \mathcal L_p\) minimization problem on randomly generated data. 
Next, we focus on the sensor network localization problem.
Finally, we apply DRSOM to the CUTEst dataset \cite{gould_cutest_2015} as a showcase in benchmark problems.

To demonstrate the efficacy of DRSOM in the experiments, we implement the algorithm using the Julia programming language for convenient comparisons to first- and second-order methods. The experiments are performed in the Julia version on a desktop of Mac OS with a 3.2 GHz 6-Core Intel Core i7 processor.
The benchmark algorithms, including the gradient descent method, conjugate gradient method, LBFGS, and Newton trust-region method, are all from a third-party package \texttt{Optim.jl}\footnote{For details, see \texttt{https://github.com/JuliaNLSolvers/Optim.jl}}, including a set of line search algorithms in \texttt{LineSearches.jl}.\footnote{For details, see \texttt{https://github.com/JuliaNLSolvers/LineSearches.jl}}. In all of our tests, we set memory parameter \cite{nocedal_numerical_2006} of LBFGS to 10.


In the experiments, we use the radius-free version of DRSOM. Let us briefly describe a strategy to update \(\mu_k\) in \eqref{eq.trs.2d.rf}. Since the ``Radius-Free'' parameter \(\mu_k\) is designed as a wild estimate to \(\lk\) for \eqref{eq.trs.2d}, we adjust \(\mu_k\) by the bounding in a prescribed box interval. Let \(\mu_1 \le\mu_2 \) be the eigenvalues of \(\Hk\) in the subspace, consider a simple adaptive rule to put \(\mu_k\) in a desired interval \([\underline{\mu}_k, \overline{\mu}_k]\):
\begin{equation}\label{eq.drsom.free.radius}
    \begin{aligned}
        \gamma_{k+1}      & = \begin{cases}
                                  \beta_2\gamma_k, \quad \rho_k \le \zeta_1 \\
                                  \max\{\underline{\gamma}, \min \left\{\sqrt{\gamma_k}, \beta_1 \gamma_k\right\}, \quad \rho_k>\zeta_2\}
                              \end{cases} \\
        \underline{\mu}_k & = \max\{0, -\mu_1\}                                                                                      \\
        \overline{\mu}_k  & = \max\{\underline{\mu}_k, \mu_2\} +\mu_M                                                                \\
        \mu_k             & =\gamma_k \cdot\overline{\mu}_k+\max \left\{1-\gamma_k, 0\right\} \cdot \underline{\mu}_k
    \end{aligned}
\end{equation}
where \(\gamma_k > 0\), \(\mu_M\) is a big number to bound \(\mu_k\) from above, and \(\underline{\gamma}\) is the minimum level for \(\gamma_k\) and that \(\beta_1 < 1 < \beta_2\). In the computations, we increase \(\gamma_k\) if the model reduction is not accurate according to \(\rho_k\).

In the above procedure \eqref{eq.drsom.free.radius}, we adjust \(\mu_k\) with \(\gamma_k\), which is expected to be less affected by the eigenvalues.
If \(\gamma_k\) approaches to \(0\), then \(\mu_k\) is close to \(\underline{\mu}_k\) which implies \(\hkap\) is almost positive semi-definite. Otherwise, \(\gamma_k\) induces a large \(\mu_k\) so as to give a small trust-region radius.

\subsection{\(\mathcal L_2 - \mathcal L_p\) Minimization}
We test the performance of DRSOM for nonconvex \(\mathcal L_2 - \mathcal L_p\)  minimization widely used in compressed sensing. Recall \(\mathcal L_2 - \mathcal L_p\)  minimization problem (\citet{ge_note_2011,chen_complexity_2014,chen_smoothing_2012}):
\begin{equation}\label{eq.app.lp}
    \min _{x \in \rea^{m}} \phi(x) = \frac{1}{2}\|A x-b\|_{2}^{2} + \lambda \|x\|_p^p,
\end{equation}
where \(A \in \rea^{n\times m}\), \(b\in \rea^n, 0 < p < 1\).
To overcome nonsmoothness of \(\|\cdot\|_p\), we apply the smoothing strategy mentioned in \citet{chen_smoothing_2012}:
\begin{equation}
    f(x)  =\frac{1}{2}\|A x-b\|_{2}^{2} + \lambda \sum_{i=1}^n s(x_i, \varepsilon)^p,
\end{equation}
where \(s(x_i, \varepsilon)\) is a piece-wise approximation of \(|x_i|\) and \(\varepsilon\) is a small pre-defined constant \(\varepsilon > 0\):
\begin{equation}
    s(x, \varepsilon)  = \begin{cases}|x| & \text { if }|x|>\varepsilon \\ \frac{x^{2}}{2 \varepsilon}+\frac{\varepsilon}{2} & \text { if }|x| \leq \varepsilon\end{cases}
\end{equation}
We randomly generate datasets with different sizes \(n, m\) based on the following procedure. The elements of matrix A are generated by \(A_{ij} \sim \mathcal N(0,1)\) with sparsity of $r=15\%, 25\%$. To construct the true sparse vector $v\in\mathbb{R}^{m}$, we let for all \(i\):
\begin{equation*}
    v_{i} \sim \begin{cases}
        0                                     & \text{with probability $p=0.5$} \\
        \mathcal N \left(0,\frac{1}{n}\right) & \text{otherwise}
    \end{cases}
\end{equation*}
Then we let $b = Av + \delta$ where $\delta$ is the noise generated as $\delta_{i} \sim \mathcal N(0,1)$ for all $i$. The parameter $\lambda$ is chosen as $\frac{1}{5}\|A^{T}b\|_{\infty}$.
We generate instances for \((n,m)\) from \((300,100)\) to \((1000, 500)\). We set \(p=0.5\) and the smoothing parameter \(\varepsilon=1e^{-1}\).

Next, we test the performance of DRSOM and benchmark algorithms including the first-order representative gradient descent method (\gd{}), and the second-order Newton trust-region method (\newtontr{}). The \gd{} is enhanced with the Hager-Zhang line-search algorithm (see \cite{zhang_nonmonotone_2004}) for global convergence. We report the running time and iteration number to reach a first-order stationary point at a precision of \(1e^{-5}\), precisely,
\begin{equation*}
    |\nabla f(\xk)| \le \epsilon := 1e^{-5}.
\end{equation*}
The iterations and running time needed for a set of methods are reported in the \autoref{tab.lp}.
\begin{table}[ht]
    \centering
    \scriptsize
    \caption{Performance of DRSOM on \eqref{eq.app.lp} compared to other algorithms on randomly generated instances: iterations needed for precision \(\epsilon = 1e^{-5}\)} \label{tab.lp}
    \begin{tabular}{r|r|r|rr|rr|rr}
        \toprule
        \multirow{2}{*}{$n$} & \multirow{2}{*}{$m$} & \multirow{2}{*}{$r$} & \multicolumn{2}{c|}{\drsom} & \multicolumn{2}{c|}{\gd} & \multicolumn{2}{c}{\newtontr}                           \\
                             &                      &                      & $k$                         & time                     & $k$                           & time    & $k$ & time    \\
        \midrule
        300                  & 100                  & 1.5e-01              & 101                         & 9.1e-02                  & 810                           & 3.8e-01 & 10  & 4.3e-02 \\
        300                  & 200                  & 1.5e-01              & 176                         & 1.3e-01                  & 3300                          & 4.3e-01 & 16  & 1.7e-01 \\
        300                  & 500                  & 1.5e-01              & 304                         & 1.8e-01                  & 3954                          & 1.2e+00 & 29  & 1.8e+00 \\
        500                  & 100                  & 1.5e-01              & 117                         & 2.3e-02                  & 818                           & 7.3e-02 & 10  & 2.0e-02 \\
        500                  & 200                  & 1.5e-01              & 199                         & 6.6e-02                  & 2682                          & 5.0e-01 & 98  & 9.7e-01 \\
        500                  & 500                  & 1.5e-01              & 306                         & 2.7e-01                  & 5301                          & 1.9e+00 & 27  & 1.8e+00 \\
        1000                 & 100                  & 1.5e-01              & 134                         & 4.1e-02                  & 985                           & 1.6e-01 & 54  & 2.0e-01 \\
        1000                 & 200                  & 1.5e-01              & 314                         & 1.4e-01                  & 2327                          & 5.9e-01 & 50  & 4.8e-01 \\
        1000                 & 500                  & 1.5e-01              & 315                         & 3.2e-01                  & 5047                          & 2.8e+00 & 57  & 3.2e+00 \\
        300                  & 100                  & 2.5e-01              & 211                         & 4.0e-02                  & 2048                          & 1.8e-01 & 70  & 1.6e-01 \\
        300                  & 200                  & 2.5e-01              & 263                         & 1.3e-01                  & 5724                          & 6.8e-01 & 69  & 7.2e-01 \\
        300                  & 500                  & 2.5e-01              & 401                         & 2.3e-01                  & 10000                         & 2.7e+00 & 125 & 6.9e+00 \\
        500                  & 100                  & 2.5e-01              & 161                         & 3.7e-02                  & 1784                          & 1.4e-01 & 55  & 1.8e-01 \\
        500                  & 200                  & 2.5e-01              & 297                         & 1.0e-01                  & 5095                          & 7.8e-01 & 90  & 9.1e-01 \\
        500                  & 500                  & 2.5e-01              & 405                         & 2.9e-01                  & 9897                          & 3.6e+00 & 209 & 1.3e+01 \\
        1000                 & 100                  & 2.5e-01              & 173                         & 1.1e-01                  & 3799                          & 5.0e-01 & 51  & 1.3e-01 \\
        1000                 & 200                  & 2.5e-01              & 286                         & 2.1e-01                  & 4954                          & 1.1e+00 & 120 & 1.3e+00 \\
        1000                 & 500                  & 2.5e-01              & 343                         & 3.6e-01                  & 7889                          & 4.3e+00 & 118 & 8.8e+00 \\
        \bottomrule
    \end{tabular}
\end{table}

These results show that the DRSOM is far better than \gd{} in all test instances in both iteration number and running time. This clearly demonstrates that the benefits of second-order information. On the other hand, compared to the full-dimensional \newtontr{}, \drsom{} achieves better running time in trade for more iterations.

\subsection{Sensor Network Localization}
We next test another nonconvex optimization problem, namely, the Sensor Network Localization (SNL). The SNL problem is to find coordinates of ad hoc wireless sensors given pairwise distances in the network. Fruitful research has been found for this problem, among which the approach based on Semidefinite Programming Relaxation (SDR) has witnessed great success; see, for example, \citet{biswas_semidefinite_2004,wang_further_2008}, and many others. We here adopt the notations in \cite{wang_further_2008}.

Let \(n\) sensors be points in \(\mathbb{R}^d\), besides, assume another set of \(m\) \emph{known} points (usually referred to as anchors) whose exact positions are \(a_1, ..., a_m\). Let \(d_{ij}\) be the distance between sensor \(i\) and \(j\), and \(\bar d_{ik}\) be the distance from the sensor \(i\) to anchor point \(k\). We can then define the set of distances as edges in the network:
\begin{equation}
    N_{x}=\left\{(i, j):\left\|x_{i}-x_{j}\right\|=d_{i j} \leq r_{d}\right\}, N_{a}=\left\{(i, k):\left\|x_{i}-a_{k}\right\|=\right. \left.d_{i k} \leq r_{d}\right\},
\end{equation}
where \(r_d\) is a fixed parameter known as the \emph{radio range}. Using these given set of distances, the SNL problem is to find other the exact positions of left $n-m$ sensors. In the following, we specifically consider $d=2$ such that the solutions to this problem can be visualized.  Formally, the SNL problem considers the following quadratic constrained quadratic programming (QCQP) feasibility problem,
\begin{subequations}
    \begin{align}
        \left\|x_{i}-x_{j}\right\|^{2}=d_{i j}^{2},       & ~ \forall(i, j) \in N_{x} \\
        \left\|x_{i}-a_{k}\right\|^{2}=\bar{d}_{i k}^{2}, & ~ \forall(i, k) \in N_{a}
    \end{align}
\end{subequations}
which can be solved by a nonlinear least-square (NLS) problem:
\begin{equation}\label{snl.nls}
    \min _{X} \sum_{(i<j, j) \in N_{x}}(\left\| x _{i}- x _{j}\right\|^{2}-d_{i j}^{2})^{2}+\sum_{(k, j) \in N_{a}}(\left\| a _{k}- x _{j}\right\|^{2}-\bar{d}_{k j}^{2})^{2},
\end{equation}
which is clearly nonconvex. Instead of directly solving \eqref{snl.nls}, one classical framework is to apply the \emph{two-stage} strategy proposed in \citet{biswas_semidefinite_2004}. In the first stage, we use semidefinite programming to solve a lifted convex relaxation to \eqref{snl.nls}:
\begin{equation}\label{snl.sdr}
    \begin{aligned}
        \mn & ~0 \bullet Z                                                                                                         \\
        \st & ~ Z_{[1:2, 1:2]}=I,                                                                                                  \\
            & \left( 0 ; e_{i}-e_{j}\right)\left( 0 ; e_{i}-e_{j}\right)^{T} \bullet Z=d_{i j}^{2} \quad \forall(i, j) \in N_{x},  \\
            & \left(-a_{k} ; e_{i}\right)\left(-a_{k} ; e_{i}\right)^{T} \bullet Z=\bar{d}_{i k}^{2} \quad \forall(i, k) \in N_{a} \\
            & Z \succeq 0,
    \end{aligned}
\end{equation}
where we let \(I\) be the identity matrix of dimension \(2\), \(e_i\) be a \(n\)-vector of zeros except for a one at \(i\)-th entry. \(Z\) is the lifted positive semidefinite matrix, such that
\begin{equation}
    Z = \begin{bmatrix}
        I   & X \\
        X^T & Y
    \end{bmatrix},
\end{equation}
which is equivalent to state \(Y \succeq X^TX\). If the solution of SDR satisfies \(\rank{Y} = 2\), then SDR solves the original problem; otherwise, \((Y,X)\) is used to initialize the NLS problem for further refinement in the second stage.  In the original work of \cite{biswas_semidefinite_2004}, local solutions of \eqref{snl.nls} are pursued by the gradient descent method (GD). To set more comparisons, we also include DRSOM, \cg{}, and \lbfgs{} in large scale problems.

\subsubsection{A specific example where \drsom{} produces better solutions}

We here provide a randomly generated example with 80 points, 5 of which are anchors. We set the radio range to 0.5 and the random distance noise \(n_f = 0.05\).  We terminate at an iterate \(\xk\) if \(\| \nabla f(\xk) \| \le 1e^{-6} \).

Our computational results of this example basically illustrate that if we initialize the NLS problem \eqref{snl.nls} for SNL by the SDR \eqref{snl.sdr}, then DRSOM and GD are comparable.  However, if we do not have the SDR solution at hand, the DRSOM may usually provide better solutions than GD, which shows the benefits of exploring second-order information in a very limited subspace.

\begin{figure}[tbhp]
    \centering
    \subfloat[The graphic result of GD initialized by an SDR solution]{
        \includegraphics[width=.43\linewidth]{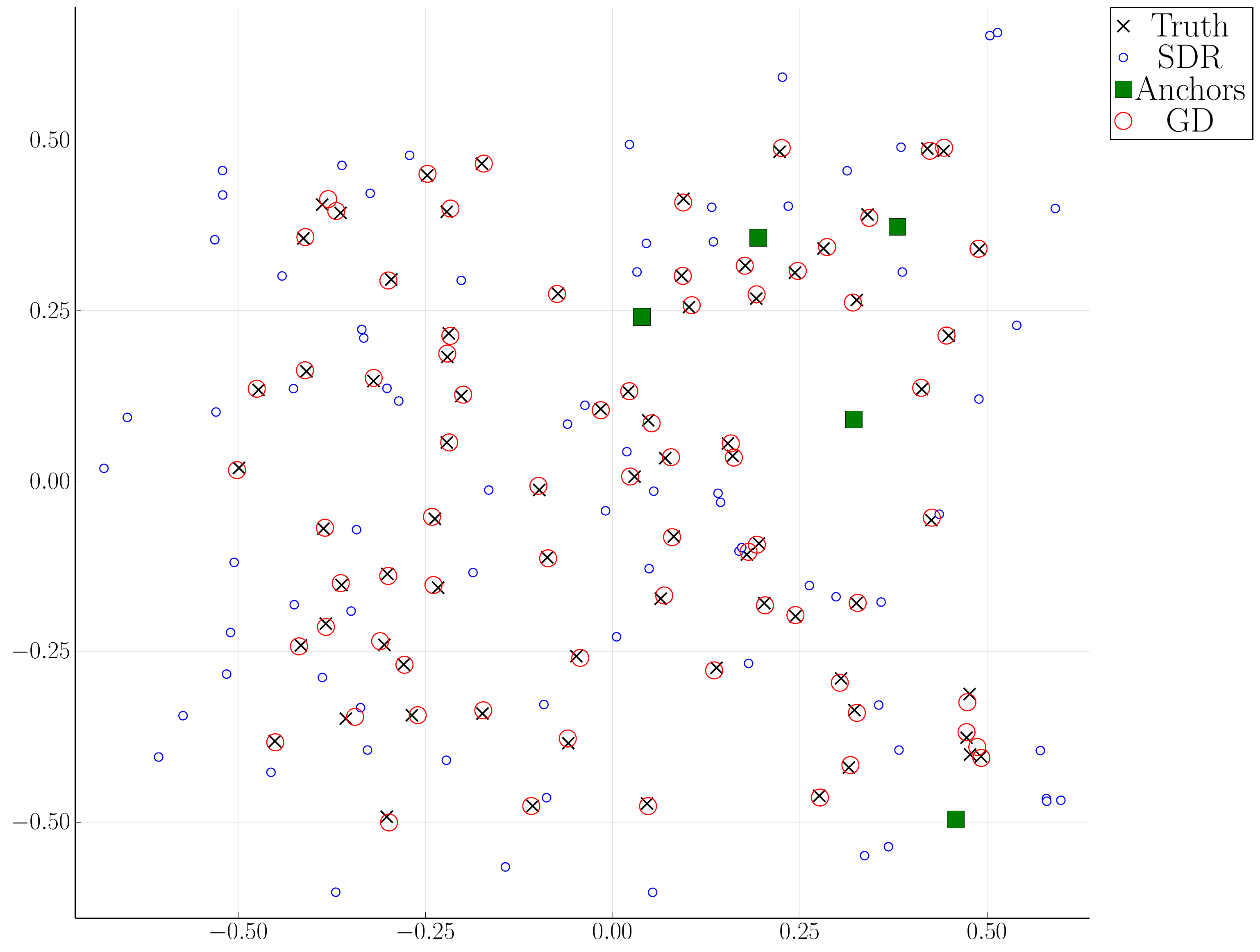}
    }\quad
    \subfloat[The graphic result of DRSOM initialized by an SDR solution]{
        \includegraphics[width=.43\linewidth]{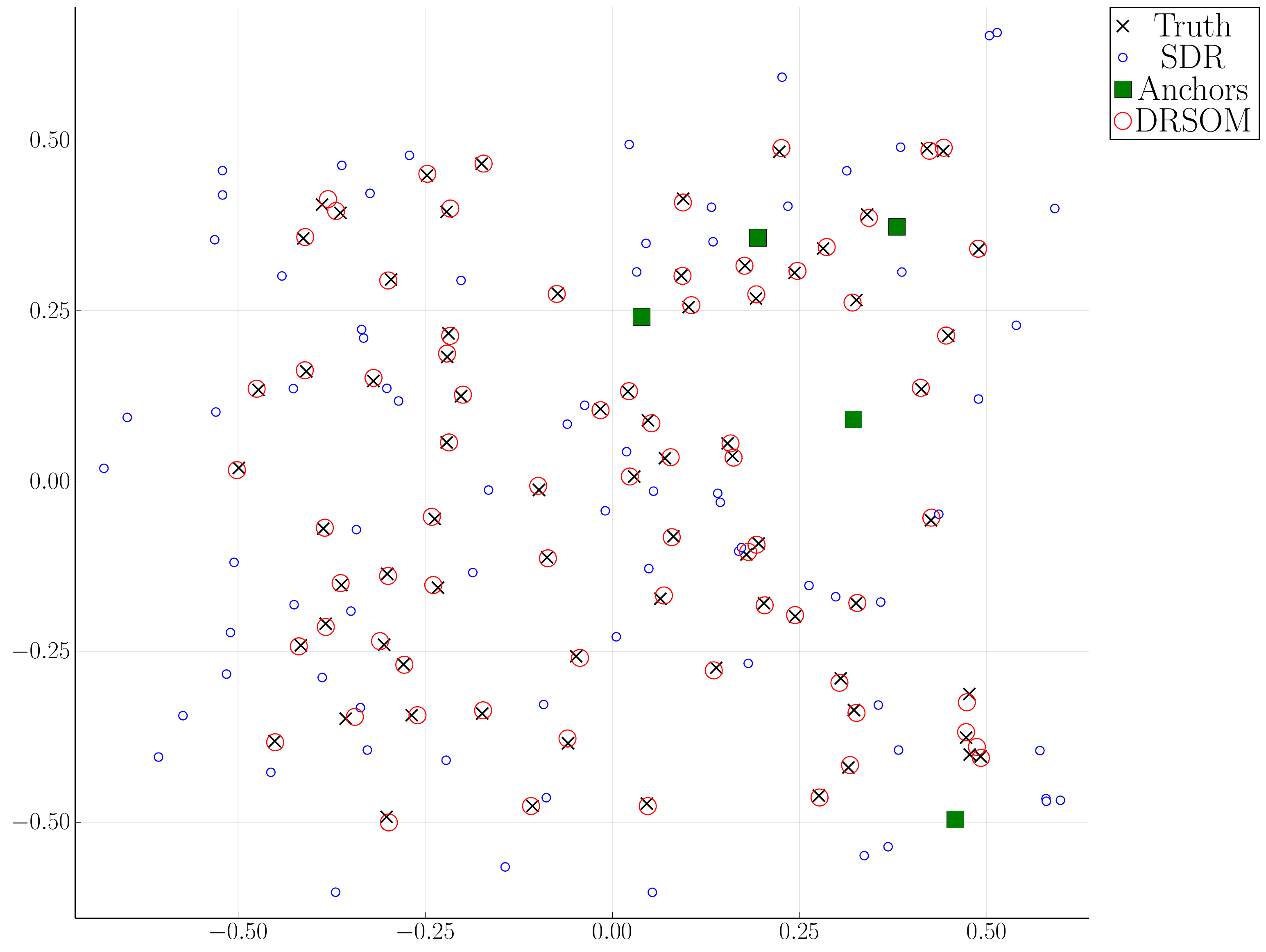}
    }
    \caption{
        Comparison of localization by GD and DRSOM with SDR initialization. The \textcolor{OliveGreen}{rectangles} and \textbf{crosses} represent the anchors and true locations, respectively. The blue \textcolor{blue}{circles} are solutions by SDR, and the red \textcolor{red}{circles} are final solutions of GD/DRSOM.
    }
    \label{fig.snl.sdr}
\end{figure}

\autoref{fig.snl.sdr} depicts the realization results of GD and DRSOM initialized by SDP relaxation. In this case, both algorithms are able to guarantee convergence to the ground truth.

\begin{figure}[h]
    \subfloat[The graphic result of GD without an SDR solution]{
        \includegraphics[width=.45\linewidth]{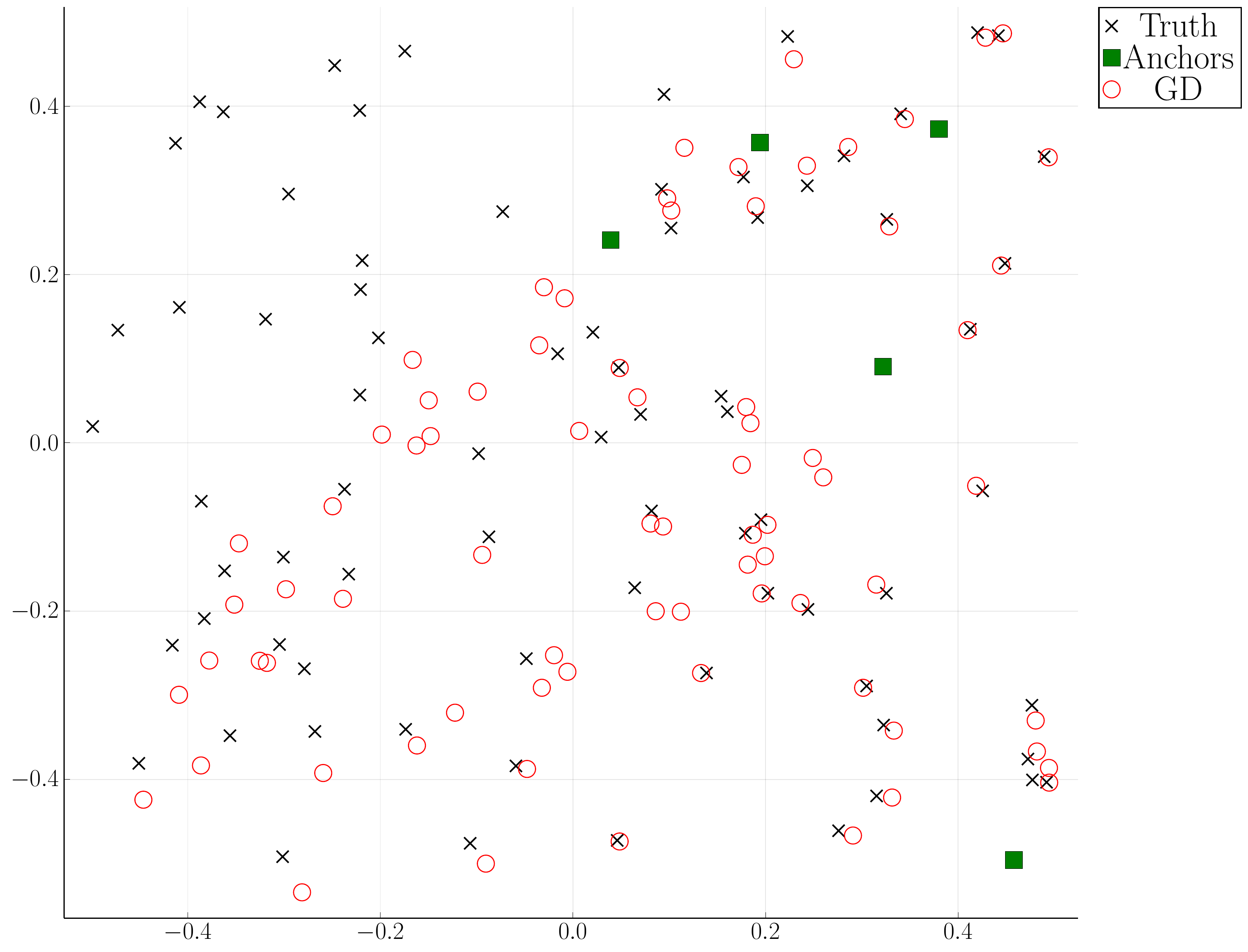}
    }\quad
    \subfloat[The graphic result of DRSOM without an SDR solution]{
        \includegraphics[width=.45\linewidth]{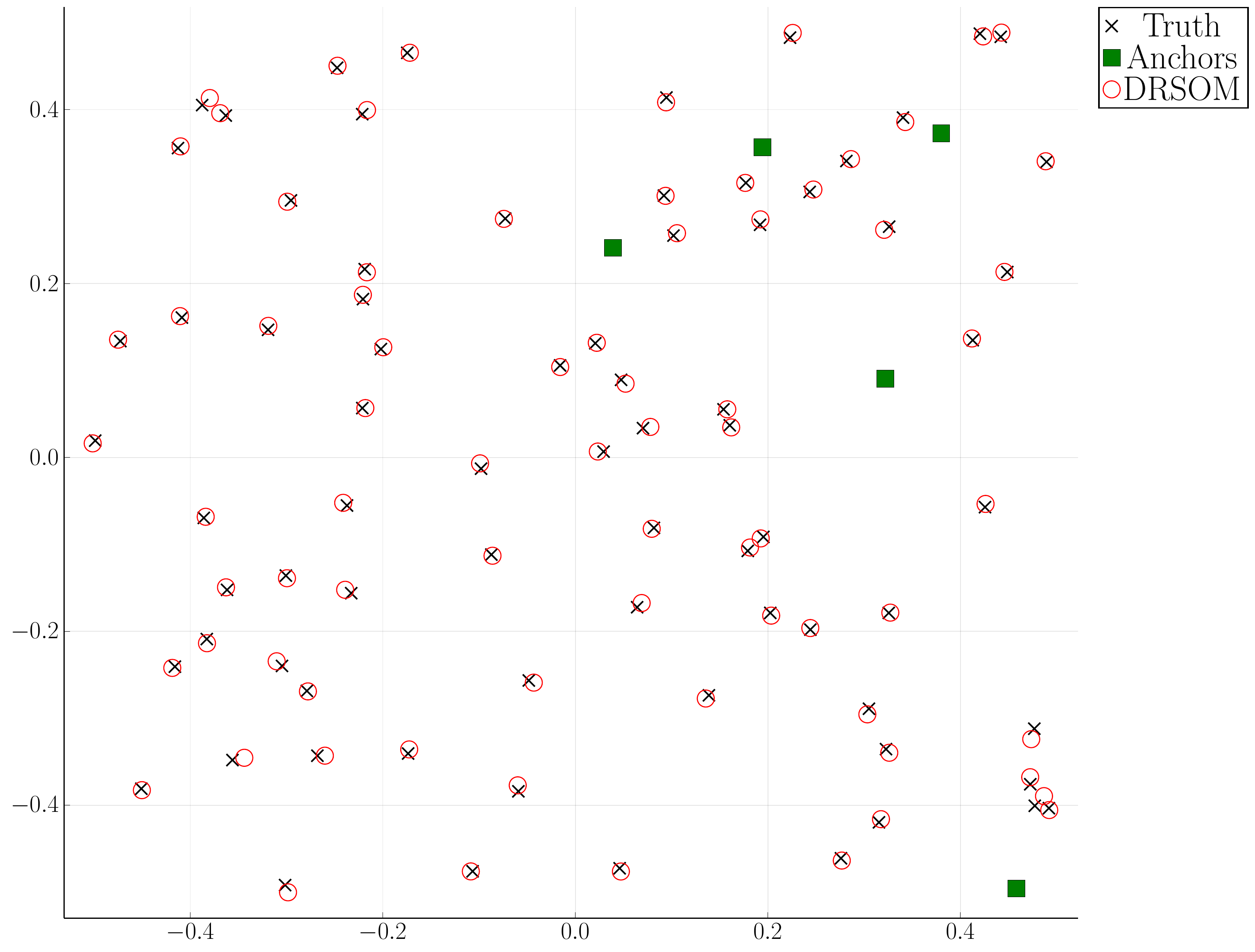}
    }
    \caption{
        Comparison of localization by GD and DRSOM without SDR initialization. The meaning of each symbol is the same as \autoref{fig.snl.sdr}.
    }
    \label{fig.snl.nosdr}
\end{figure}

As a comparison, \autoref{fig.snl.nosdr} depicts the case without solving SDR first. We use the same parameter settings as for the previous case. The GD and DRSOM initializes with \(X_{i}:= 0, i=1,...,n\) and both solves to \(\|\nabla f(\xk) \| \le 1e^{-6} \).
In this very particular case, the GD fails to recover true positions; however, the DRSOM can provide accurate solutions even without the initialization of SDR.

In this example, we rigorously provide a case where the DRSOM may result in a better solution (in this case, the global one) than the first-order methods despite that we only have optimal subspace guarantees in theory.

\subsubsection{Solving large SNL instances}
In this part, we compare \gd{}, \cg{}, \lbfgs{} and \drsom{} on large SNL instances up to 10,000 sensors.  Similar to the previous tests, \cg{}, \lbfgs{} and \gd{} are facilitated with Hager-Zhang line search algorithm. We terminate an algorithm if $\|\nabla f(\xk)\| \le 1e^{-5}$ or the running time exceeds 3,000 seconds.

Our results are presented in \autoref{tab:snl-large.k}. The \gd{} fails to converge if the number of sensors exceeds 3,000, so we do not try it on larger problems. In all instances, \drsom{} stands out in running time and has a clear advantage over \cg{} and \lbfgs{}. Notably, it successfully solves the problem with 10,000 sensors to $\|\nabla f(\xk)\|\le 1e^{-5}$ in about 2,200 seconds while all other competing algorithms fail.

\begin{table}[]
    \footnotesize
    \caption{Results of \cg{}, \lbfgs{}, \drsom{}, and \gd{} on SNL instances. $n$ and $m$ denotes the number of sensors and anchors. $|E|$ denotes the number of used edges (QCQP constraints).} \label{tab:snl-large.k}
    \centering
    \subfloat[Iteration number]{
        \begin{tabular}{ccc|cccc}
            \toprule
            \multirow{2}{*}{$n$} & \multirow{2}{*}{$m$} & \multirow{2}{*}{$|E|$} & \multicolumn{4}{c}{$k$}                                \\
                                 &                      &                        & \cg                     & \lbfgs{} & \drsom  & \gd     \\
            \midrule
            500                  & 50                   & 2.2e+04                & 1.2e+02                 & 9.5e+01  & 9.0e+01 & 3.2e+02 \\
            1000                 & 80                   & 4.6e+04                & 2.3e+02                 & 2.1e+02  & 2.1e+02 & 1.4e+03 \\
            2000                 & 120                  & 9.4e+04                & 4.5e+02                 & 3.8e+02  & 3.6e+02 & 4.4e+03 \\
            3000                 & 150                  & 1.4e+05                & 7.4e+02                 & 5.6e+02  & 5.4e+02 & 7.5e+03 \\
            4000                 & 400                  & 1.8e+05                & 1.2e+03                 & 5.7e+02  & 7.8e+02 & -       \\
            6000                 & 600                  & 2.7e+05                & 1.2e+03                 & 1.1e+03  & 1.1e+03 & -       \\
            10000                & 1000                 & 4.5e+05                & 1.0e+03                 & 9.7e+02  & 1.5e+03 & -       \\
            \bottomrule
        \end{tabular}
    }

    \subfloat[Gradient norm at termination. Symbol ``*'' means the algorithm does not attain the required accuracy.]{
        \begin{tabular}{ccc|cccc}
            \toprule
            \multirow{2}{*}{$n$} & \multirow{2}{*}{$m$} & \multirow{2}{*}{$|E|$} & \multicolumn{4}{c}{$\|g\|$}                                       \\
                                 &                      &                        & \cg                         & \lbfgs{}    & \drsom  & \gd         \\
            \midrule
            500                  & 50                   & 2.2e+04                & 7.0e-06                     & 8.9e-06     & 6.8e-06 & 9.4e-06     \\
            1000                 & 80                   & 4.6e+04                & 7.0e-06                     & 8.3e-06     & 8.6e-06 & 9.9e-06     \\
            2000                 & 120                  & 9.4e+04                & 8.4e-06                     & 8.5e-06     & 8.4e-06 & 9.9e-06     \\
            3000                 & 150                  & 1.4e+05                & 8.9e-06                     & 9.3e-06     & 8.4e-06 & 1.6e-03$^*$ \\
            4000                 & 400                  & 1.8e+05                & 7.9e-06                     & 8.9e-06     & 9.8e-06 & -           \\
            6000                 & 600                  & 2.7e+05                & 9.2e-06                     & 9.9e-06     & 8.8e-06 & -           \\
            10000                & 1000                 & 4.5e+05                & 2.3e-03$^*$                 & 1.0e-04$^*$ & 9.0e-06 & -           \\
            \bottomrule
        \end{tabular}
    }

    \subfloat[Running time. Symbol ``-'' means the algorithm reaches 3,000s]{
        \begin{tabular}{ccc|cccc}
            \toprule
            \multirow{2}{*}{$n$} & \multirow{2}{*}{$m$} & \multirow{2}{*}{$|E|$} & \multicolumn{4}{c}{$t$}                                         \\
                                 &                      &                        & \cg{}                   & \lbfgs{} & \drsom{}         & \gd{}   \\
            \midrule
            500                  & 50                   & 2.2e+04                & 1.7e+01                 & 1.3e+01  & \textbf{1.1e+01} & 2.3e+01 \\
            1000                 & 80                   & 4.6e+04                & 7.3e+01                 & 6.2e+01  & \textbf{3.9e+01} & 1.8e+02 \\
            2000                 & 120                  & 9.4e+04                & 2.5e+02                 & 2.5e+02  & \textbf{1.4e+02} & 1.1e+03 \\
            3000                 & 150                  & 1.4e+05                & 6.5e+02                 & 5.8e+02  & \textbf{2.7e+02} & -       \\
            4000                 & 400                  & 1.8e+05                & 1.3e+03                 & 7.2e+02  & \textbf{5.0e+02} & -       \\
            6000                 & 600                  & 2.7e+05                & 2.0e+03                 & 2.1e+03  & \textbf{1.1e+03} & -       \\
            10000                & 1000                 & 4.5e+05                & -                       & -        & \textbf{2.2e+03} & -       \\
            \bottomrule
        \end{tabular}
    }
    \normalsize
\end{table}

\subsection{CUTEst Benchmark}
In this subsection, we test DRSOM in a nonlinear programming benchmark CUTEst dataset \cite{gould_cutest_2015}. Besides, we report results of first-order methods including \gd{} and \cg{}, and second-order methods including \newtontr{} and \lbfgs{}. Similar to the previous tests, \cg{}, \lbfgs{} and \gd{} are facilitated with Hager-Zhang line search algorithm. Since each problem in the CUTEst dataset may have multiple different realizations by possibly multiple parameters, we enumerate these combinations and choose the first one that has $n \le 200$ variables.

We use the nonlinear programming interface in Julia programming language facilitated by open source infrastructure \cite{orban_juliasmoothoptimizers_2019} for DRSOM and all other competing algorithms. We report an overall comparison in \autoref{tab.cutest.overall}, where we present average values of iteration number $\overline k$, function $\overline k^f$, gradient $\overline k^g$ and Hessian $\overline k^H$ evaluations.
Besides, we also provide the \emph{scaled geometric mean} as an extra metric including $\overline k_G$, $\overline k_G^f$, $\overline k_G^g$, $\overline k_G^H$ correspondingly. We report CPU time $\overline t, \overline t_G$ in both measures. We mark an instance successfully solved if the returned iterate $\xk$ satisfies $\min\{\|\gk\|,\|\gk\|/\|g_0\|\} \le 1e^{-5}$, where $g_0$ is the gradient at the initial point $x_0$. The minimum we use here accounts for the case where the gradient is too large, in which case a solution is acceptable in a relative measure. Under such a criterion, we count the total number of successful instances as $\mathcal K$.

Note that the Hessian-vector products are already available in \cite{orban_juliasmoothoptimizers_2019}, so DRSOM directly takes advantage of HVP computations. To this end, the metric $\overline k_G^g$ also includes number of needed HVP as extra gradient evaluations (counted as 2 per each call). We scale geometric mean for time and all other evaluations number by 1 (second) and 50 (iterations), respectively. We report detailed results in \autoref{sec.cutest.appendix}.

\begin{table}[ht]
    \centering
    \caption{Performance of different algorithms on the CUTEst dataset.
        Note $\overline t, \overline k$ are mean running time and iteration;
        $\overline t_{G}, \overline k_{G}$ are scaled geometric means (scaled by 1 second and 50 iterations, respectively).
        If an instance is the failed, its iteration number and solving time are set to $20,000$. }\label{tab.cutest.overall}
    \subfloat[Performance in algebraic averages]{
        \scriptsize
        \begin{tabular}{lrrrrrr}
            \toprule
            method    & $\mathcal K$ & $\overline t$ & $\overline k$ & $\overline k^f$ & $\overline k^g$ & $\overline k^H$ \\
            \midrule
            \lbfgs    & 96.00        & 0.18          & 583.53        & 1680.79         & 1680.79         & 0.00            \\
            \newtontr & 94.00        & 1.38          & 2636.71       & 927.78          & 2637.71         & 875.96          \\
            \cg       & 98.00        & 0.45          & 1193.41       & 12374.70        & 12374.70        & 0.00            \\
            \drsom    & 97.00        & 0.20          & 1163.92       & 1206.06         & 3519.81         & 0.00            \\
            \gd       & 91.00        & 1.33          & 6007.38       & 24602.65        & 24602.65        & 0.00            \\
            \bottomrule
        \end{tabular}
    }

    \subfloat[Performance in geometric averages]{
        \scriptsize
        \begin{tabular}{lrrrrrr}
            \toprule
            method    & $\overline t_G$ & $\overline k_G$ & $\overline k_G^f$ & $\overline k_G^g$ & $\overline k_G^H$ \\
            \midrule
            \lbfgs    & 0.08            & 76.39           & 180.59            & 180.59            & -0.00             \\
            \newtontr & 0.29            & 99.34           & 43.84             & 101.30            & 40.35             \\
            \cg       & 0.09            & 124.73          & 282.38            & 282.38            & -0.00             \\
            \drsom    & 0.09            & 133.10          & 162.42            & 332.32            & -0.00             \\
            \gd       & 0.32            & 952.60          & 1947.67           & 1947.67           & -0.00             \\
            \bottomrule
        \end{tabular}
    }
\end{table}

From \autoref{tab.cutest.overall}, we can see that the quasi-Newton method \lbfgs{} stands out in terms of iteration number, while \cg{} has the most successful instances. However, we observe that \drsom{} is exceptional using only 2D subspace with gradient and momentum. This observation becomes more evident in average CPU time $\overline t$ and evaluation numbers.
Not surprisingly, \autoref{tab.cutest.overall} indicates that DRSOM has a clear advantage over the gradient method, since \drsom{} solves more instances and tends to be more stable, as seen from scaled metrics with subscript $G$. Besides, the average running time of \drsom{} is better than \cg{}, but these two methods should be close in terms of scaled metrics. We highlight a comparison of \drsom{} with \lbfgs{} (with a Hager-Zhang line-search algorithm), which clearly demonstrates the benefits of using second-order information even in a very limited manner. As the current implementation of \drsom{} already demonstrates a saving in function and gradient evaluations, it is reasonable to expect better numerical routines can further improve these benchmark results. This of course serves as one of our future directions.
\section{Conclusion}
In this paper, we introduce a Dimension-Reduced Second-Order Method (DRSOM) that utilizes 2-D trust-region subproblems and requires only the gradient, momentum and two extra Hessian-vector products (alternatively, interpolation with three function evaluations) that cost almost the same as first-order methods.
We give a concise analysis of the global and local speed of convergence.  Computationally, we provide fruitful examples in nonconvex optimization. Notably, DRSOM has substantial superiority to first-order methods while keeping a low cost per iteration compared to second-order methods.

Our preliminary results motivate future developments of DRSOM. In terms of \autoref{assm.hessian}, it is an important aspect to find more efficient approach to satisfy \eqref{eq.assm.hessian} in addition to the Krylov-like corrector method \autoref{alg.corrector}.

\section*{Acknowledgement}
\addcontentsline{toc}{section}{Acknowledgements}
The authors would like to thank Tianyi Lin for the insightful comments that significantly improved the presentation of this paper.

\addcontentsline{toc}{section}{References}
\clearpage

\clearpage
\bibliography{main}
\bibliographystyle{plainnat-ex}
\clearpage
\appendix
\begin{landscape}
  \section{Detailed Computational Results for CUTEst Dataset}\label{sec.cutest.appendix}
  \scriptsize
  \begin{longtable}{llllllllllll}
\caption{Complete Results on CUTEst Dataset, iteration \& time}
\label{tab.cutest.kt}\\
\toprule
        & method & \multicolumn{2}{l}{\cg} & \multicolumn{2}{l}{\drsom} & \multicolumn{2}{l}{\gd} & \multicolumn{2}{l}{\lbfgs} & \multicolumn{2}{l}{\newtontr} \\
        & {} &    $k$ &      $t$ &    $k$ &      $t$ &    $k$ &      $t$ &    $k$ &      $t$ &       $k$ &      $t$ \\
name & n &        &          &        &          &        &          &        &          &           &          \\
\midrule
\endfirsthead
\caption[]{Complete Results on CUTEst Dataset, iteration \& time} \\
\toprule
        & method & \multicolumn{2}{l}{\cg} & \multicolumn{2}{l}{\drsom} & \multicolumn{2}{l}{\gd} & \multicolumn{2}{l}{\lbfgs} & \multicolumn{2}{l}{\newtontr} \\
        & {} &    $k$ &      $t$ &    $k$ &      $t$ &    $k$ &      $t$ &    $k$ &      $t$ &       $k$ &      $t$ \\
name & n &        &          &        &          &        &          &        &          &           &          \\
\midrule
\endhead
\midrule
\multicolumn{12}{r}{{Continued on next page}} \\
\midrule
\endfoot

\bottomrule
\endlastfoot
ARGLINA & 200 &      1 &  4.5e-04 &      2 &  1.0e-03 &      1 &  4.4e-04 &      1 &  6.1e-04 &         5 &  4.9e-01 \\
ARGLINB & 200 &      1 &  4.7e-04 &      2 &  2.0e-03 &      1 &  5.7e-04 &      3 &  1.5e-03 &         4 &  3.8e-01 \\
ARGLINC & 200 &      1 &  4.4e-04 &      2 &  2.0e-03 &      1 &  5.1e-04 &      4 &  2.4e-03 &         4 &  3.4e-01 \\
ARGTRIGLS & 200 &   1151 &  8.7e-01 &    611 &  1.0e+00 &  20000 &  1.7e+01 &    541 &  6.4e-01 &        12 &  1.2e+00 \\
ARWHEAD & 100 &     18 &  8.2e-04 &     10 &  3.0e-03 &     57 &  2.6e-03 &      6 &  5.6e-04 &         5 &  5.1e-03 \\
BDQRTIC & 100 &    136 &  6.4e-03 &    138 &  2.1e-02 &   3701 &  1.5e-01 &     30 &  1.5e-03 &        10 &  1.8e-02 \\
BOX & 10  &      7 &  1.9e-04 &      4 &  0.0e+00 &     60 &  6.3e-04 &      4 &  2.7e-04 &         2 &  3.6e-04 \\
BOXPOWER & 10  &     59 &  6.1e-04 &     19 &  1.0e-03 &  20000 &  1.3e-01 &     15 &  6.0e-04 &        11 &  5.9e-04 \\
BROWNAL & 200 &      5 &  8.9e-03 &      5 &  1.1e-02 &    102 &  6.6e-02 &      4 &  5.9e-03 &         6 &  5.6e-01 \\
BROYDN3DLS & 50  &     25 &  4.4e-04 &     30 &  3.0e-03 &     62 &  8.3e-04 &     24 &  7.4e-04 &         5 &  2.3e-03 \\
BROYDN7D & 50  &     85 &  3.0e-03 &     82 &  9.0e-03 &    167 &  4.9e-03 &     63 &  2.5e-03 &        15 &  6.8e-03 \\
BROYDNBDLS & 50  &    113 &  2.5e-03 &    102 &  1.1e-02 &    996 &  2.5e-02 &     61 &  2.0e-03 &        10 &  4.8e-03 \\
BRYBND & 50  &    113 &  2.9e-03 &    102 &  1.4e-02 &    996 &  2.1e-02 &     61 &  2.0e-03 &        10 &  4.4e-03 \\
CHAINWOO & 4   &    201 &  1.4e-03 &    326 &  1.6e-02 &   6053 &  2.7e-02 &     21 &  3.2e-04 &        45 &  8.5e-04 \\
CHNROSNB & 25  &    352 &  3.8e-03 &    322 &  2.0e-02 &   4754 &  5.8e-02 &    109 &  1.6e-03 &        40 &  4.2e-03 \\
CHNRSNBM & 25  &    361 &  4.0e-03 &    374 &  2.4e-02 &   6252 &  6.5e-02 &    171 &  2.4e-03 &        45 &  4.8e-03 \\
COSINE & 100 &      2 &  1.8e-04 &      2 &  0.0e+00 &      2 &  4.6e-04 &     10 &  9.3e-04 &        12 &  2.8e-02 \\
CRAGGLVY & 50  &     96 &  3.4e-03 &     87 &  1.3e-02 &    330 &  1.4e-02 &     71 &  2.6e-03 &        17 &  6.3e-03 \\
CURLY10 & 100 &   1274 &  4.1e-02 &    820 &  5.5e-02 &  17334 &  5.4e-01 &    809 &  2.5e-02 &        12 &  2.1e-02 \\
CURLY20 & 100 &   1306 &  5.6e-02 &    933 &  6.7e-02 &  20000 &  8.9e-01 &    882 &  3.4e-02 &        12 &  2.5e-02 \\
DIXMAANA & 90  &      7 &  3.6e-04 &      8 &  2.0e-03 &      7 &  2.8e-04 &      5 &  5.4e-04 &         7 &  1.7e-02 \\
DIXMAANB & 90  &      8 &  4.1e-04 &      9 &  2.0e-03 &      9 &  3.5e-04 &      6 &  6.1e-04 &        12 &  7.5e-02 \\
DIXMAANC & 90  &      9 &  4.7e-04 &     10 &  2.0e-03 &      9 &  5.2e-04 &      6 &  5.0e-04 &        10 &  1.5e-02 \\
DIXMAAND & 90  &     10 &  5.2e-04 &     10 &  2.0e-03 &     11 &  7.1e-04 &      8 &  5.8e-04 &        12 &  1.8e-02 \\
DIXMAANE & 90  &     43 &  2.5e-03 &     47 &  7.0e-03 &    364 &  1.6e-02 &     42 &  1.7e-03 &         8 &  1.2e-02 \\
DIXMAANF & 90  &     43 &  2.5e-03 &     49 &  8.0e-03 &    362 &  1.7e-02 &     37 &  1.6e-03 &        15 &  2.0e-02 \\
DIXMAANG & 90  &     44 &  2.6e-03 &     49 &  8.0e-03 &    324 &  1.3e-02 &     35 &  1.6e-03 &        13 &  2.0e-02 \\
DIXMAANH & 90  &     42 &  2.5e-03 &     49 &  8.0e-03 &    355 &  1.4e-02 &     37 &  1.6e-03 &        15 &  2.6e-02 \\
DIXMAANI & 90  &    151 &  8.8e-03 &    155 &  2.4e-02 &  14646 &  6.5e-01 &    154 &  5.9e-03 &        10 &  1.3e-02 \\
DIXMAANJ & 90  &    136 &  8.0e-03 &    142 &  2.0e-02 &  12973 &  6.9e-01 &    138 &  5.5e-03 &        17 &  3.0e-02 \\
DIXMAANK & 90  &    137 &  8.2e-03 &    134 &  1.6e-02 &  13414 &  6.0e-01 &    138 &  5.4e-03 &        16 &  3.3e-02 \\
DIXMAANL & 90  &    121 &  7.1e-03 &    129 &  1.7e-02 &  13357 &  6.9e-01 &    138 &  5.3e-03 &        17 &  2.2e-02 \\
DIXMAANM & 90  &    161 &  1.9e-02 &    162 &  1.2e+00 &  15424 &  7.4e-01 &    160 &  8.6e-01 &         8 &  2.8e-01 \\
DIXMAANN & 90  &    202 &  1.2e-02 &    202 &  2.5e-02 &  14490 &  7.1e-01 &    167 &  6.4e-03 &        20 &  2.9e-02 \\
DIXMAANO & 90  &    147 &  8.6e-03 &    202 &  2.4e-02 &  14655 &  7.6e-01 &    181 &  7.0e-03 &        19 &  5.7e-02 \\
DIXMAANP & 90  &    149 &  8.9e-03 &    188 &  2.4e-02 &  14575 &  7.6e-01 &    136 &  5.3e-03 &        26 &  3.3e-02 \\
DIXON3DQ & 100 &    100 &  1.9e-03 &    102 &  1.1e-02 &  20000 &  3.1e-01 &    100 &  2.2e-03 &         5 &  9.1e-03 \\
DQDRTIC & 50  &      5 &  1.9e-04 &      9 &  1.0e-03 &   1333 &  2.8e-02 &      5 &  4.6e-04 &         5 &  1.5e-03 \\
DQRTIC & 50  &     14 &  3.1e-04 &     16 &  1.0e-03 &     14 &  3.3e-04 &     31 &  6.4e-04 &        25 &  4.7e-03 \\
EDENSCH & 36  &     30 &  8.0e-04 &     27 &  3.0e-03 &     44 &  1.6e-03 &     21 &  6.6e-04 &        15 &  4.1e-03 \\
EIGENALS & 6   &     17 &  3.7e-04 &      7 &  1.0e-03 &    349 &  1.9e-03 &      9 &  3.6e-04 &         7 &  4.6e-04 \\
EIGENBLS & 6   &     32 &  2.6e-04 &     30 &  2.0e-03 &    282 &  3.0e-03 &     11 &  3.1e-04 &        17 &  1.2e-03 \\
EIGENCLS & 30  &    175 &  3.9e-03 &    133 &  1.4e-02 &    841 &  2.0e-02 &    126 &  3.4e-03 &        17 &  4.4e-03 \\
ENGVAL1 & 50  &     24 &  5.8e-04 &     24 &  3.0e-03 &     36 &  5.9e-04 &     15 &  5.1e-04 &         9 &  1.1e-02 \\
ERRINROS & 25  &  19275 &  1.9e-01 &   6025 &  2.3e-01 &  20000 &  2.5e-01 &     94 &  1.5e-03 &        47 &  4.9e-03 \\
ERRINRSM & 25  &  20000 &  1.9e-01 &  20000 &  5.8e-01 &  20000 &  2.1e-01 &    217 &  2.9e-03 &        86 &  8.9e-03 \\
EXTROSNB & 100 &   3798 &  8.6e-02 &   2722 &  1.4e-01 &  20000 &  5.4e-01 &   5037 &  1.3e-01 &      1645 &  1.7e+00 \\
FLETBV3M & 10  &      0 &  8.1e-06 &      3 &  0.0e+00 &      0 &  9.5e-07 &      0 &  5.2e-05 &         0 &  5.2e-05 \\
FLETCBV2 & 10  &      5 &  1.6e-04 &      6 &  1.0e-03 &    182 &  2.0e-03 &      5 &  2.7e-04 &         1 &  3.3e-04 \\
FLETCBV3 & 10  &      0 &  9.1e-06 &     14 &  2.0e-03 &      0 &  6.0e-06 &      0 &  8.5e-05 &         0 &  4.8e-05 \\
FLETCHBV & 10  &     10 &  1.1e-04 &     20 &  2.0e-03 &    390 &  2.9e-03 &     10 &  3.2e-04 &        88 &  2.9e-03 \\
FLETCHCR & 100 &    678 &  3.3e-02 &    866 &  1.1e-01 &  20000 &  5.4e-01 &    487 &  1.5e-02 &       202 &  5.0e-01 \\
FMINSRF2 & 16  &     68 &  6.1e-04 &     57 &  4.0e-03 &    983 &  7.5e-03 &     29 &  5.8e-04 &     20000 &  7.0e-02 \\
FMINSURF & 16  &     24 &  2.1e-04 &     27 &  3.0e-03 &     84 &  8.6e-04 &     21 &  4.2e-04 &     20000 &  1.7e+00 \\
FREUROTH & 50  &     97 &  2.6e-03 &     99 &  1.2e-02 &   6575 &  1.6e-01 &     16 &  1.0e-03 &         9 &  5.2e-03 \\
GENHUMPS & 10  &    362 &  6.8e-03 &    134 &  1.2e-02 &   5166 &  9.6e-02 &    421 &  6.0e-03 &     12910 &  6.1e-01 \\
GENROSE & 100 &    309 &  1.7e-02 &    272 &  2.8e-02 &   2924 &  9.8e-02 &    270 &  7.3e-02 &        95 &  2.0e-01 \\
HILBERTA & 6   &      4 &  1.9e-04 &      6 &  1.0e-03 &   2839 &  2.4e-02 &      4 &  2.7e-04 &         4 &  4.9e-04 \\
HILBERTB & 5   &      3 &  6.6e-05 &      4 &  0.0e+00 &      6 &  4.9e-05 &      3 &  3.4e-04 &         3 &  3.4e-04 \\
INDEF & 50  &      1 &  9.1e-06 &     31 &  3.8e-02 &      1 &  3.1e-06 &      1 &  8.5e-05 &     20000 &  9.3e+00 \\
INDEFM & 50  &    135 &  4.1e-03 &    152 &  1.5e-02 &  20000 &  7.4e-01 &     63 &  2.4e-03 &        28 &  1.5e-02 \\
INTEQNELS & 102 &      5 &  1.3e-03 &      6 &  4.0e-03 &      6 &  1.2e-03 &      4 &  1.3e-03 &         3 &  3.4e-02 \\
JIMACK & 81  &  12150 &  3.4e+01 &   3285 &  1.0e+01 &  20000 &  7.1e+01 &   3628 &  1.5e+01 &     20000 &  5.1e+01 \\
LIARWHD & 36  &     24 &  6.6e-04 &     12 &  1.0e-03 &    337 &  6.0e-03 &     10 &  3.9e-04 &        12 &  4.0e-03 \\
MANCINO & 50  &      9 &  1.2e-02 &     10 &  2.4e-02 &     10 &  1.4e-02 &      8 &  1.2e-02 &        11 &  3.3e-02 \\
MODBEALE & 10  &    194 &  2.1e-03 &   2222 &  7.5e-02 &  20000 &  2.8e-01 &     36 &  7.4e-04 &         8 &  5.5e-04 \\
MOREBV & 50  &    218 &  3.5e-03 &    220 &  1.9e-02 &  20000 &  4.3e-01 &    252 &  5.3e-03 &         1 &  5.6e-04 \\
MSQRTALS & 49  &    135 &  4.6e-03 &    132 &  2.0e-02 &   1184 &  4.5e-02 &    103 &  5.0e-03 &        12 &  1.0e-02 \\
MSQRTBLS & 49  &    217 &  7.2e-03 &    249 &  3.0e-02 &   3570 &  1.6e-01 &    177 &  7.8e-03 &        15 &  1.3e-02 \\
NCB20 & 110 &    852 &  1.0e-01 &    731 &  1.9e-01 &  20000 &  6.8e+00 &    330 &  4.8e-02 &        75 &  1.9e-01 \\
NCB20B & 180 &   2336 &  4.8e-01 &   4710 &  2.2e+00 &  20000 &  1.1e+01 &   1154 &  3.5e-01 &        11 &  5.9e-02 \\
NONCVXU2 & 10  &     50 &  4.0e-04 &     51 &  5.0e-03 &    182 &  1.9e-03 &     23 &  4.5e-04 &     20000 &  1.6e+00 \\
NONCVXUN & 10  &     29 &  5.4e-04 &     28 &  2.0e-03 &     58 &  5.0e-04 &     21 &  4.3e-04 &     20000 &  8.1e-01 \\
NONDIA & 90  &     19 &  6.9e-04 &     10 &  2.0e-03 &   7974 &  1.7e-01 &      8 &  4.9e-04 &         6 &  4.2e-03 \\
NONDQUAR & 100 &   3040 &  7.2e-02 &   7558 &  3.2e-01 &   8735 &  2.3e-01 &    666 &  1.7e-02 &        15 &  1.2e-01 \\
NONMSQRT & 49  &  20000 &  5.7e-01 &  20000 &  1.3e+00 &  20000 &  9.4e-01 &  20000 &  9.7e-01 &     20000 &  8.6e+00 \\
OSCIGRAD & 15  &     73 &  1.9e-03 &     63 &  5.0e-03 &    183 &  3.2e-03 &     45 &  9.4e-04 &     20000 &  3.2e+00 \\
OSCIPATH & 25  &      9 &  2.4e-04 &     10 &  2.0e-03 &     12 &  1.1e-04 &      9 &  4.5e-04 &         9 &  2.1e-03 \\
PENALTY1 & 50  &     44 &  8.9e-04 &     37 &  4.0e-03 &    673 &  1.7e-02 &     90 &  1.5e-03 &        41 &  1.9e-02 \\
PENALTY2 & 50  &    194 &  8.1e-03 &    172 &  2.2e-02 &   5765 &  2.2e-01 &    139 &  5.9e-03 &        25 &  1.3e-02 \\
PENALTY3 & 50  &     84 &  1.6e-01 &     90 &  2.5e-01 &    329 &  6.0e-01 &     53 &  1.0e-01 &        17 &  6.0e-02 \\
POWELLSG & 60  &     80 &  1.1e-03 &   1017 &  4.6e-02 &  11094 &  1.4e-01 &     23 &  6.0e-04 &        16 &  5.2e-03 \\
POWER & 50  &     34 &  4.8e-04 &     31 &  2.0e-03 &     53 &  1.3e-03 &     27 &  5.7e-04 &        21 &  1.0e-02 \\
QUARTC & 100 &     16 &  5.4e-04 &     20 &  3.0e-03 &     16 &  4.2e-04 &     31 &  8.1e-04 &        29 &  1.6e-02 \\
SBRYBND & 50  &  20000 &  9.5e+00 &  20000 &  1.0e+00 &  20000 &  1.2e+01 &  20000 &  7.4e-01 &       576 &  3.3e-01 \\
SCHMVETT & 10  &     32 &  6.9e-04 &     28 &  2.0e-03 &    211 &  2.3e-03 &     17 &  3.9e-04 &         3 &  3.9e-04 \\
SCOSINE & 10  &      2 &  1.5e-04 &      2 &  1.0e-03 &      2 &  1.3e-04 &      2 &  3.5e-04 &     20000 &  6.9e-01 \\
SCURLY10 & 10  &      1 &  1.0e-05 &      2 &  1.3e-02 &      1 &  3.1e-06 &      1 &  5.4e-05 &        11 &  9.4e-04 \\
SENSORS & 10  &     18 &  5.0e-03 &     19 &  5.0e-03 &     28 &  4.1e-03 &     25 &  2.3e-03 &        10 &  1.1e-03 \\
SINQUAD & 50  &     60 &  2.0e-03 &     25 &  4.0e-03 &    444 &  2.0e-02 &     12 &  7.3e-04 &        10 &  4.2e-03 \\
SPARSINE & 50  &    275 &  8.8e-03 &    320 &  2.9e-02 &   7976 &  2.2e-01 &    196 &  6.2e-03 &        28 &  1.3e-02 \\
SPARSQUR & 50  &     10 &  2.6e-04 &     13 &  3.0e-03 &     10 &  4.0e-04 &     11 &  6.7e-04 &        14 &  8.5e-03 \\
SPMSRTLS & 100 &     65 &  3.5e-03 &     65 &  1.3e-02 &    286 &  1.4e-02 &     60 &  3.4e-03 &        13 &  2.0e-02 \\
SROSENBR & 50  &     14 &  4.5e-04 &      8 &  1.0e-03 &   5273 &  5.2e-02 &      7 &  3.3e-04 &         9 &  2.0e-03 \\
SSBRYBND & 50  &  11679 &  6.2e-01 &   3751 &  1.8e-01 &  20000 &  3.5e+00 &   2528 &  9.2e-02 &        21 &  9.4e-03 \\
SSCOSINE & 10  &      8 &  1.0e-03 &  20000 &  4.9e-01 &      2 &  1.2e-04 &      2 &  4.2e-04 &     20000 &  5.4e-01 \\
TOINTGSS & 50  &     16 &  3.6e-04 &     17 &  3.0e-03 &     37 &  1.3e-03 &     16 &  5.8e-04 &     20000 &  1.3e+01 \\
TQUARTIC & 50  &     25 &  5.6e-04 &     20 &  3.0e-03 &   9119 &  1.6e-01 &     11 &  5.1e-04 &        12 &  4.8e-03 \\
TRIDIA & 50  &     52 &  6.0e-04 &     59 &  4.0e-03 &   3187 &  3.2e-02 &     52 &  1.0e-03 &         4 &  1.6e-03 \\
VARDIM & 200 &     10 &  7.5e-03 &      8 &  1.1e-02 &      4 &  2.7e-03 &      3 &  8.3e-04 &         9 &  4.8e-02 \\
VAREIGVL & 100 &     98 &  5.4e-03 &    126 &  2.1e-02 &     55 &  2.3e-03 &     21 &  1.3e-03 &        25 &  5.5e-02 \\
WATSON & 12  &    681 &  1.2e-02 &    386 &  3.1e-02 &  20000 &  5.3e-01 &     38 &  1.1e-03 &        13 &  1.9e-03 \\
WOODS & 4   &    201 &  1.5e-03 &    326 &  1.6e-02 &   6053 &  3.8e-02 &     21 &  3.4e-04 &        45 &  8.3e-04 \\
YATP1LS & 120 &    357 &  1.4e-01 &     58 &  2.6e-02 &  20000 &  2.8e+00 &    129 &  1.9e-02 &     20000 &  4.6e+01 \\
YATP2LS & 8   &      9 &  1.2e-04 &     11 &  2.0e-03 &     13 &  1.6e-04 &      8 &  2.9e-04 &     20000 &  5.1e-01 \\
\end{longtable}

\begin{longtable}{llllllllllll}
\caption{Complete Results on CUTEst Dataset, function value \& norm of the gradient}
\label{tab.cutest.fx}\\
\toprule
        & method & \multicolumn{2}{l}{\cg} & \multicolumn{2}{l}{\drsom} & \multicolumn{2}{l}{\gd} & \multicolumn{2}{l}{\lbfgs} & \multicolumn{2}{l}{\newtontr} \\
        & {} &  $\|g\|$ &       $f$ &  $\|g\|$ &       $f$ &  $\|g\|$ &       $f$ &  $\|g\|$ &       $f$ &   $\|g\|$ &       $f$ \\
name & n &          &           &          &           &          &           &          &           &           &           \\
\midrule
\endfirsthead
\caption[]{Complete Results on CUTEst Dataset, function value \& norm of the gradient} \\
\toprule
        & method & \multicolumn{2}{l}{\cg} & \multicolumn{2}{l}{\drsom} & \multicolumn{2}{l}{\gd} & \multicolumn{2}{l}{\lbfgs} & \multicolumn{2}{l}{\newtontr} \\
        & {} &  $\|g\|$ &       $f$ &  $\|g\|$ &       $f$ &  $\|g\|$ &       $f$ &  $\|g\|$ &       $f$ &   $\|g\|$ &       $f$ \\
name & n &          &           &          &           &          &           &          &           &           &           \\
\midrule
\endhead
\midrule
\multicolumn{12}{r}{{Continued on next page}} \\
\midrule
\endfoot

\bottomrule
\endlastfoot
ARGLINA & 200 &  4.0e-14 &   2.3e-26 &  3.0e-13 &   2.3e-26 &  4.0e-14 &   2.3e-26 &  4.0e-14 &   2.3e-26 &   5.2e-14 &   1.3e-26 \\
ARGLINB & 200 &  2.8e-03 &   5.0e+01 &  2.3e-02 &   5.0e+01 &  2.8e-03 &   5.0e+01 &  7.6e-06 &   5.0e+01 &   1.9e-03 &   5.0e+01 \\
ARGLINC & 200 &  1.5e-03 &   5.1e+01 &  1.3e-02 &   5.1e+01 &  1.5e-03 &   5.1e+01 &  7.1e-06 &   5.1e+01 &   1.5e-03 &   5.1e+01 \\
ARGTRIGLS & 200 &  9.8e-06 &   2.6e-12 &  9.4e-06 &   1.7e-14 &  4.6e-05 &   5.2e-10 &  9.0e-06 &   1.4e-12 &   4.4e-09 &   8.3e-20 \\
ARWHEAD & 100 &  8.5e-07 &   3.0e-12 &  1.2e-07 &   0.0e+00 &  7.5e-06 &   2.1e-11 &  1.4e-07 &   0.0e+00 &   6.3e-06 &   6.6e-14 \\
BDQRTIC & 100 &  7.2e-06 &   3.8e+02 &  9.4e-06 &   3.8e+02 &  1.0e-05 &   3.8e+02 &  8.0e-06 &   3.8e+02 &   4.0e-06 &   3.8e+02 \\
BOX & 10  &  4.7e-07 &  -1.7e-01 &  5.6e-06 &  -1.7e-01 &  7.3e-06 &  -1.7e-01 &  1.8e-08 &  -1.7e-01 &   2.6e-06 &  -1.7e-01 \\
BOXPOWER & 10  &  8.6e-06 &   1.2e-07 &  6.5e-06 &   4.3e-09 &  7.0e-05 &   3.5e-06 &  4.1e-08 &   6.0e-11 &   5.6e-06 &   1.3e-07 \\
BROWNAL & 200 &  4.8e-06 &   1.5e-09 &  3.8e-06 &   1.5e-09 &  9.3e-06 &   1.5e-09 &  1.6e-06 &   1.5e-09 &   4.7e-10 &   7.5e-20 \\
BROYDN3DLS & 50  &  6.9e-06 &   2.4e-12 &  9.3e-06 &   5.2e-13 &  8.0e-06 &   1.8e-11 &  5.9e-06 &   2.2e-12 &   1.0e-06 &   5.7e-14 \\
BROYDN7D & 50  &  7.6e-06 &   1.8e+01 &  9.7e-06 &   1.8e+01 &  9.8e-06 &   1.7e+01 &  6.8e-06 &   1.7e+01 &   1.9e-09 &   1.7e+01 \\
BROYDNBDLS & 50  &  1.0e-05 &   4.1e-12 &  8.7e-06 &   2.1e-13 &  9.9e-06 &   2.2e-11 &  1.0e-05 &   1.1e-11 &   4.4e-13 &   8.8e-18 \\
BRYBND & 50  &  1.0e-05 &   4.1e-12 &  8.7e-06 &   2.1e-13 &  9.9e-06 &   2.2e-11 &  1.0e-05 &   1.1e-11 &   4.4e-13 &   8.8e-18 \\
CHAINWOO & 4   &  3.5e-06 &   1.0e+00 &  9.0e-06 &   1.0e+00 &  9.9e-06 &   1.0e+00 &  2.9e-07 &   1.0e+00 &   2.7e-09 &   1.0e+00 \\
CHNROSNB & 25  &  7.1e-06 &   1.9e-11 &  7.3e-06 &   3.1e-12 &  1.0e-05 &   1.4e-10 &  3.7e-06 &   4.2e-13 &   4.6e-07 &   2.5e-13 \\
CHNRSNBM & 25  &  1.0e-05 &   7.3e-12 &  9.0e-06 &   3.5e-12 &  1.0e-05 &   1.3e-10 &  7.8e-06 &   2.6e-12 &   2.7e-09 &   2.0e-20 \\
COSINE & 100 &  8.7e+00 &  -1.3e+01 &  8.1e+01 &  -1.3e+01 &  8.7e+00 &  -1.3e+01 &  2.1e-06 &  -9.9e+01 &   2.1e-13 &  -9.9e+01 \\
CRAGGLVY & 50  &  9.6e-06 &   1.5e+01 &  8.4e-06 &   1.5e+01 &  9.9e-06 &   1.5e+01 &  9.7e-06 &   1.5e+01 &   1.6e-07 &   1.5e+01 \\
CURLY10 & 100 &  2.4e-05 &  -1.0e+04 &  1.1e-04 &  -1.0e+04 &  3.3e-05 &  -1.0e+04 &  1.8e-05 &  -1.0e+04 &   3.3e-10 &  -1.0e+04 \\
CURLY20 & 100 &  2.6e-05 &  -1.0e+04 &  1.5e-04 &  -1.0e+04 &  7.4e-05 &  -1.0e+04 &  3.3e-05 &  -1.0e+04 &   2.2e-10 &  -1.0e+04 \\
DIXMAANA & 90  &  3.0e-07 &   1.0e+00 &  8.7e-07 &   1.0e+00 &  2.0e-06 &   1.0e+00 &  4.8e-06 &   1.0e+00 &   6.0e-19 &   1.0e+00 \\
DIXMAANB & 90  &  3.5e-06 &   1.0e+00 &  9.8e-07 &   1.0e+00 &  3.3e-07 &   1.0e+00 &  3.1e-06 &   1.0e+00 &   2.5e-07 &   1.0e+00 \\
DIXMAANC & 90  &  2.8e-07 &   1.0e+00 &  3.6e-07 &   1.0e+00 &  1.2e-06 &   1.0e+00 &  5.2e-06 &   1.0e+00 &   5.1e-12 &   1.0e+00 \\
DIXMAAND & 90  &  2.5e-06 &   1.0e+00 &  9.0e-06 &   1.0e+00 &  4.3e-06 &   1.0e+00 &  6.9e-07 &   1.0e+00 &   4.3e-11 &   1.0e+00 \\
DIXMAANE & 90  &  6.7e-06 &   1.0e+00 &  6.7e-06 &   1.0e+00 &  9.7e-06 &   1.0e+00 &  8.3e-06 &   1.0e+00 &   1.8e-12 &   1.0e+00 \\
DIXMAANF & 90  &  7.6e-06 &   1.0e+00 &  8.6e-06 &   1.0e+00 &  9.8e-06 &   1.0e+00 &  9.7e-06 &   1.0e+00 &   2.7e-07 &   1.0e+00 \\
DIXMAANG & 90  &  7.0e-06 &   1.0e+00 &  7.5e-06 &   1.0e+00 &  9.9e-06 &   1.0e+00 &  6.6e-06 &   1.0e+00 &   1.0e-11 &   1.0e+00 \\
DIXMAANH & 90  &  9.9e-06 &   1.0e+00 &  8.8e-06 &   1.0e+00 &  9.7e-06 &   1.0e+00 &  9.2e-06 &   1.0e+00 &   2.7e-10 &   1.0e+00 \\
DIXMAANI & 90  &  8.9e-06 &   1.0e+00 &  9.4e-06 &   1.0e+00 &  1.0e-05 &   1.0e+00 &  9.6e-06 &   1.0e+00 &   3.2e-09 &   1.0e+00 \\
DIXMAANJ & 90  &  8.7e-06 &   1.0e+00 &  8.5e-06 &   1.0e+00 &  1.0e-05 &   1.0e+00 &  6.7e-06 &   1.0e+00 &   1.2e-06 &   1.0e+00 \\
DIXMAANK & 90  &  1.0e-05 &   1.0e+00 &  9.7e-06 &   1.0e+00 &  1.0e-05 &   1.0e+00 &  6.7e-06 &   1.0e+00 &   6.5e-07 &   1.0e+00 \\
DIXMAANL & 90  &  8.9e-06 &   1.0e+00 &  9.4e-06 &   1.0e+00 &  1.0e-05 &   1.0e+00 &  9.8e-06 &   1.0e+00 &   1.0e-11 &   1.0e+00 \\
DIXMAANM & 90  &  9.4e-06 &   1.0e+00 &  9.2e-06 &   1.0e+00 &  1.0e-05 &   1.0e+00 &  8.7e-06 &   1.0e+00 &   2.0e-15 &   1.0e+00 \\
DIXMAANN & 90  &  8.3e-06 &   1.0e+00 &  9.5e-06 &   1.0e+00 &  1.0e-05 &   1.0e+00 &  9.9e-06 &   1.0e+00 &   5.7e-07 &   1.0e+00 \\
DIXMAANO & 90  &  9.6e-06 &   1.0e+00 &  9.8e-06 &   1.0e+00 &  1.0e-05 &   1.0e+00 &  8.7e-06 &   1.0e+00 &   9.0e-06 &   1.0e+00 \\
DIXMAANP & 90  &  8.1e-06 &   1.0e+00 &  9.8e-06 &   1.0e+00 &  1.0e-05 &   1.0e+00 &  8.9e-06 &   1.0e+00 &   2.2e-07 &   1.0e+00 \\
DIXON3DQ & 100 &  6.3e-12 &   3.9e-23 &  8.4e-06 &   2.1e-12 &  1.5e-04 &   4.6e-04 &  5.5e-12 &   2.6e-23 &   2.5e-14 &   4.9e-25 \\
DQDRTIC & 50  &  3.2e-14 &   6.4e-29 &  4.4e-06 &   3.5e-12 &  9.9e-06 &   2.5e-11 &  4.1e-13 &   7.2e-27 &   4.5e-14 &   1.2e-28 \\
DQRTIC & 50  &  1.6e-06 &   1.3e-07 &  4.8e-06 &   3.4e-09 &  9.7e-06 &   9.8e-07 &  2.8e-06 &   2.3e-08 &   8.9e-06 &   4.3e-07 \\
EDENSCH & 36  &  5.8e-06 &   2.2e+02 &  3.0e-06 &   2.2e+02 &  9.9e-06 &   2.2e+02 &  8.7e-06 &   2.2e+02 &   6.0e-06 &   2.2e+02 \\
EIGENALS & 6   &  1.1e-07 &   1.6e-16 &  1.7e-06 &   9.4e-24 &  9.8e-06 &   9.1e-11 &  9.9e-07 &   1.3e-14 &   1.4e-06 &   5.6e-13 \\
EIGENBLS & 6   &  4.5e-06 &   1.8e-01 &  8.6e-06 &   1.8e-01 &  9.9e-06 &   1.8e-01 &  3.3e-08 &   1.8e-01 &   5.6e-09 &   4.9e-18 \\
EIGENCLS & 30  &  9.6e-06 &   1.2e-10 &  7.0e-06 &   5.6e-12 &  9.9e-06 &   6.1e-10 &  5.5e-06 &   4.5e-11 &   1.2e-07 &   4.8e-16 \\
ENGVAL1 & 50  &  6.5e-06 &   5.4e+01 &  7.7e-06 &   5.4e+01 &  9.6e-06 &   5.4e+01 &  4.3e-06 &   5.4e+01 &   1.4e-08 &   5.4e+01 \\
ERRINROS & 25  &  9.2e-06 &   1.8e+01 &  9.8e-06 &   1.8e+01 &  4.0e-03 &   1.8e+01 &  9.9e-06 &   1.8e+01 &   6.9e-06 &   1.8e+01 \\
ERRINRSM & 25  &  1.5e-03 &   1.8e+01 &  1.1e-01 &   1.8e+01 &  1.4e-02 &   1.8e+01 &  9.7e-06 &   1.8e+01 &   3.7e-07 &   1.8e+01 \\
EXTROSNB & 100 &  9.9e-06 &   3.1e-06 &  8.9e-06 &   1.9e-06 &  6.8e-03 &   1.6e-03 &  8.3e-06 &   2.8e-09 &   9.7e-06 &   1.8e-08 \\
FLETBV3M & 10  &  2.4e-06 &   1.2e-05 &  5.6e-07 &  -2.2e-03 &  2.4e-06 &   1.2e-05 &  2.4e-06 &   1.2e-05 &   2.4e-06 &   1.2e-05 \\
FLETCBV2 & 10  &  6.1e-09 &  -5.5e-01 &  3.6e-07 &  -5.5e-01 &  1.0e-05 &  -5.5e-01 &  6.1e-09 &  -5.5e-01 &   1.7e-15 &  -5.5e-01 \\
FLETCBV3 & 10  &  2.4e-06 &   1.2e-05 &  3.5e-07 &  -3.2e-02 &  2.4e-06 &   1.2e-05 &  2.4e-06 &   1.2e-05 &   2.4e-06 &   1.2e-05 \\
FLETCHBV & 10  &  9.1e-13 &  -2.7e+06 &  2.7e-06 &  -2.7e+06 &  9.3e-06 &  -2.7e+06 &  0.0e+00 &  -2.7e+06 &   4.5e-13 &  -2.7e+06 \\
FLETCHCR & 100 &  9.3e-06 &   3.2e-11 &  9.2e-06 &   1.7e-12 &  1.0e-02 &   1.4e-04 &  6.9e-06 &   9.9e-13 &   2.5e-06 &   2.7e-12 \\
FMINSRF2 & 16  &  7.0e-06 &   1.0e+00 &  9.0e-06 &   1.0e+00 &  1.0e-05 &   1.0e+00 &  5.2e-06 &   1.0e+00 &   6.5e-01 &   0.0e+00 \\
FMINSURF & 16  &  1.0e-05 &   1.0e+00 &  6.1e-06 &   1.0e+00 &  9.4e-06 &   1.0e+00 &  6.4e-06 &   1.0e+00 &   8.8e-01 &   3.6e+01 \\
FREUROTH & 50  &  4.6e-06 &   5.9e+03 &  6.3e-06 &   5.9e+03 &  3.9e-05 &   5.9e+03 &  6.0e-06 &   5.9e+03 &   7.1e-08 &   5.9e+03 \\
GENHUMPS & 10  &  3.5e-06 &   9.6e-11 &  2.4e+01 &   4.4e+02 &  4.1e-06 &   9.6e-11 &  9.6e-06 &   4.5e-10 &   2.8e-07 &   4.7e-13 \\
GENROSE & 100 &  8.0e-06 &   1.0e+00 &  5.7e-06 &   1.0e+00 &  7.2e-06 &   1.0e+00 &  6.6e-06 &   1.0e+00 &   3.8e-07 &   1.0e+00 \\
HILBERTA & 6   &  2.2e-07 &   5.0e-09 &  1.2e-06 &   5.0e-09 &  1.0e-05 &   1.4e-07 &  2.2e-07 &   5.0e-09 &   2.6e-16 &   4.0e-25 \\
HILBERTB & 5   &  2.2e-06 &   7.5e-13 &  3.9e-06 &   2.6e-17 &  8.6e-07 &   6.5e-14 &  2.2e-06 &   7.5e-13 &   2.3e-14 &   6.2e-29 \\
INDEF & 50  &  1.8e+00 &   4.6e+01 &  8.6e+00 &  -9.3e+14 &  1.8e+00 &   4.6e+01 &  1.8e+00 &   4.6e+01 &   1.1e+00 &  -7.2e+15 \\
INDEFM & 50  &  7.3e-06 &  -4.8e+03 &  8.4e-06 &  -4.6e+03 &  1.1e-04 &  -4.6e+03 &  6.7e-07 &  -4.9e+03 &   6.2e-06 &  -5.0e+03 \\
INTEQNELS & 102 &  2.5e-06 &   4.3e-11 &  1.7e-06 &   1.8e-15 &  2.6e-06 &   4.7e-11 &  6.1e-06 &   1.8e-10 &   2.5e-11 &   3.0e-21 \\
JIMACK & 81  &  9.1e-06 &   8.8e-01 &  1.0e-05 &   8.9e-01 &  1.8e-02 &   8.7e-01 &  8.8e-06 &   9.2e-01 &   1.3e+01 &   1.2e+00 \\
LIARWHD & 36  &  7.1e-07 &   1.0e-15 &  7.0e-08 &   7.1e-30 &  9.7e-06 &   3.3e-11 &  1.8e-10 &   2.0e-19 &   2.5e-09 &   8.2e-20 \\
MANCINO & 50  &  3.4e-06 &   8.2e-17 &  8.2e-07 &   5.2e-22 &  7.3e-07 &   3.1e-18 &  4.5e-08 &   5.2e-21 &   1.6e-09 &   6.5e-24 \\
MODBEALE & 10  &  6.7e-06 &   4.7e-11 &  9.8e-06 &   1.2e-10 &  4.2e-04 &   2.7e-07 &  2.5e-07 &   4.8e-16 &   1.4e-07 &   7.9e-15 \\
MOREBV & 50  &  4.2e-06 &   8.2e-10 &  7.2e-06 &   7.6e-10 &  1.0e-05 &   6.0e-06 &  8.3e-06 &   1.5e-08 &   5.8e-06 &   6.4e-09 \\
MSQRTALS & 49  &  9.6e-06 &   1.1e-10 &  9.3e-06 &   1.3e-11 &  9.9e-06 &   3.4e-09 &  1.0e-05 &   5.1e-10 &   6.0e-08 &   1.3e-14 \\
MSQRTBLS & 49  &  8.5e-06 &   2.1e-10 &  9.7e-06 &   4.2e-11 &  1.0e-05 &   1.3e-08 &  7.0e-06 &   2.7e-10 &   5.2e-08 &   9.9e-15 \\
NCB20 & 110 &  6.6e-06 &   1.8e+02 &  1.0e-05 &   1.8e+02 &  3.9e-04 &   1.9e+02 &  9.4e-06 &   1.8e+02 &   4.1e-06 &   1.8e+02 \\
NCB20B & 180 &  8.6e-06 &   3.5e+02 &  1.0e-05 &   3.5e+02 &  1.2e-04 &   3.5e+02 &  9.6e-06 &   3.5e+02 &   4.7e-06 &   3.5e+02 \\
NONCVXU2 & 10  &  4.4e-06 &   2.3e+01 &  9.5e-06 &   2.3e+01 &  9.6e-06 &   2.3e+01 &  3.4e-06 &   2.3e+01 &   7.7e-01 &   2.3e+01 \\
NONCVXUN & 10  &  3.4e-06 &   2.3e+01 &  9.8e-06 &   2.3e+01 &  6.8e-06 &   2.3e+01 &  3.7e-06 &   2.3e+01 &   7.2e-04 &   2.6e+01 \\
NONDIA & 90  &  3.9e-07 &   3.9e-15 &  7.3e-09 &   5.0e-28 &  1.0e-05 &   9.0e-12 &  3.1e-08 &   2.8e-18 &   1.5e-07 &   4.1e-18 \\
NONDQUAR & 100 &  9.6e-06 &   5.6e-07 &  1.0e-05 &   1.6e-08 &  1.0e-05 &   4.0e-05 &  1.0e-05 &   1.9e-06 &   4.7e-06 &   2.7e-09 \\
NONMSQRT & 49  &  1.9e-02 &   1.1e+00 &  7.8e-02 &   1.1e+00 &  1.0e-02 &   1.2e+00 &  2.8e-03 &   1.1e+00 &   1.0e-01 &   1.1e+00 \\
OSCIGRAD & 15  &  5.3e-06 &   2.8e-09 &  8.7e-06 &   2.8e-09 &  9.6e-06 &   2.8e-09 &  3.4e-06 &   2.8e-09 &   2.3e-04 &   1.2e-09 \\
OSCIPATH & 25  &  4.7e-06 &   1.0e+00 &  3.8e-06 &   1.0e+00 &  4.6e-06 &   1.0e+00 &  3.8e-06 &   1.0e+00 &   1.3e-12 &   1.0e+00 \\
PENALTY1 & 50  &  6.3e-06 &   4.3e-04 &  1.6e-06 &   4.3e-04 &  8.2e-06 &   4.3e-04 &  2.6e-06 &   4.3e-04 &   8.7e-06 &   4.3e-04 \\
PENALTY2 & 50  &  9.3e-06 &   4.3e+00 &  3.1e-06 &   4.3e+00 &  1.0e-05 &   4.3e+00 &  6.5e-06 &   4.3e+00 &   2.2e-10 &   4.3e+00 \\
PENALTY3 & 50  &  9.2e-06 &   1.0e-03 &  8.0e-06 &   1.0e-03 &  9.1e-06 &   1.0e-03 &  8.6e-06 &   1.0e-03 &   7.5e-07 &   1.0e-03 \\
POWELLSG & 60  &  6.1e-06 &   2.5e-07 &  9.9e-06 &   1.4e-09 &  1.0e-05 &   7.1e-07 &  8.1e-06 &   9.0e-11 &   3.2e-06 &   5.7e-08 \\
POWER & 50  &  9.0e-06 &   3.1e-08 &  6.0e-06 &   2.3e-09 &  7.1e-06 &   2.1e-08 &  7.8e-06 &   4.6e-09 &   3.8e-06 &   6.7e-09 \\
QUARTC & 100 &  9.7e-06 &   3.1e-06 &  8.5e-06 &   1.2e-08 &  9.8e-06 &   2.8e-06 &  6.6e-06 &   6.2e-08 &   4.9e-06 &   4.4e-07 \\
SBRYBND & 50  &  7.2e+05 &   6.5e+02 &  1.6e+04 &   2.8e+02 &  8.0e+04 &   7.4e+02 &  4.2e+05 &   4.5e+02 &   4.1e-06 &   3.7e-13 \\
SCHMVETT & 10  &  7.9e-06 &  -2.4e+01 &  9.2e-06 &  -2.4e+01 &  9.7e-06 &  -2.4e+01 &  3.8e-07 &  -2.4e+01 &   1.8e-13 &  -2.4e+01 \\
SCOSINE & 10  &  2.0e+14 &   6.4e-01 &  2.0e+14 &   6.4e-01 &  3.4e+13 &   7.8e-01 &  3.4e+13 &   7.8e-01 &   5.3e+04 &  -6.8e+00 \\
SCURLY10 & 10  &  1.1e+26 &   4.2e+26 &  4.4e+97 &   4.2e+26 &  1.1e+26 &   4.2e+26 &  1.1e+26 &   4.2e+26 &   1.0e+21 &   8.4e+19 \\
SENSORS & 10  &  4.3e-06 &  -2.1e+01 &  7.8e-06 &  -2.1e+01 &  3.0e-06 &  -2.1e+01 &  1.3e-06 &  -2.1e+01 &   4.0e-08 &  -2.0e+01 \\
SINQUAD & 50  &  7.7e-06 &  -1.1e+03 &  9.4e-06 &  -2.0e+02 &  9.9e-06 &  -1.1e+03 &  3.1e-08 &  -1.1e+03 &   8.3e-10 &  -1.1e+03 \\
SPARSINE & 50  &  7.4e-06 &   9.2e-12 &  9.2e-06 &   1.5e-12 &  1.0e-05 &   4.1e-11 &  8.8e-06 &   1.1e-11 &   5.9e-10 &   1.4e-19 \\
SPARSQUR & 50  &  9.4e-06 &   1.1e-07 &  5.6e-07 &   9.2e-11 &  3.0e-06 &   2.2e-08 &  9.1e-06 &   1.1e-07 &   6.9e-06 &   6.3e-08 \\
SPMSRTLS & 100 &  6.9e-06 &   3.3e-10 &  9.1e-06 &   1.5e-11 &  9.9e-06 &   1.6e-09 &  8.5e-06 &   3.0e-10 &   5.3e-07 &   4.8e-14 \\
SROSENBR & 50  &  6.9e-06 &   9.1e-12 &  5.3e-11 &   3.4e-21 &  1.0e-05 &   2.0e-09 &  3.3e-06 &   2.6e-10 &   9.7e-09 &   2.1e-18 \\
SSBRYBND & 50  &  8.6e-06 &   2.1e-13 &  9.6e-06 &   1.2e-14 &  1.8e+02 &   4.7e+01 &  9.1e-06 &   8.9e-14 &   1.8e-09 &   1.1e-17 \\
SSCOSINE & 10  &  1.1e+11 &  -7.0e+00 &  8.1e+02 &  -8.9e+00 &  1.5e+07 &  -1.1e+00 &  1.5e+07 &  -1.1e+00 &   3.5e+02 &  -6.9e+00 \\
TOINTGSS & 50  &  9.2e-06 &   1.0e+01 &  8.0e-06 &   1.0e+01 &  8.6e-06 &   1.0e+01 &  5.0e-06 &   1.0e+01 &   1.5e+00 &   1.1e+01 \\
TQUARTIC & 50  &  8.5e-07 &   1.8e-13 &  6.9e-08 &   5.9e-14 &  1.0e-05 &   6.5e-09 &  2.9e-06 &   1.2e-14 &   1.7e-12 &   1.2e-25 \\
TRIDIA & 50  &  3.2e-06 &   1.5e-13 &  7.8e-06 &   7.8e-14 &  1.0e-05 &   5.1e-11 &  3.2e-06 &   1.5e-13 &   1.7e-14 &   1.4e-29 \\
VARDIM & 200 &  2.3e+10 &   8.8e+09 &  6.9e+10 &   8.9e+08 &  1.4e+10 &   4.4e+09 &  4.7e+10 &   2.3e+10 &   7.8e+10 &   4.5e+10 \\
VAREIGVL & 100 &  9.9e-06 &   1.9e-09 &  9.8e-06 &   1.6e-10 &  7.7e-06 &   1.2e-10 &  7.3e-06 &   5.5e-11 &   7.7e-09 &   1.2e-17 \\
WATSON & 12  &  7.8e-06 &   2.3e-07 &  7.3e-06 &   1.6e-07 &  2.7e-04 &   1.1e-05 &  1.1e-06 &   1.6e-07 &   1.9e-06 &   3.2e-07 \\
WOODS & 4   &  3.5e-06 &   1.4e-12 &  9.0e-06 &   3.5e-12 &  9.9e-06 &   1.7e-10 &  2.9e-07 &   1.4e-16 &   2.7e-09 &   1.2e-19 \\
YATP1LS & 120 &  2.9e-06 &   3.8e-08 &  9.0e-06 &   6.7e-18 &  6.2e-03 &   4.3e-02 &  1.3e-07 &   1.7e-17 &   9.8e+01 &   8.0e+01 \\
YATP2LS & 8   &  2.9e-06 &   7.5e-13 &  6.1e-07 &   3.0e-16 &  2.4e-06 &   3.6e-13 &  1.1e-09 &   9.4e-20 &   4.0e-02 &   1.0e-04 \\
\end{longtable}

  \normalsize
\end{landscape}

\end{document}